\documentclass[12pt]{amsart}
\usepackage{amssymb, amsmath, amsfonts, amsthm, graphics, mathrsfs}
\usepackage[hmargin=1 in, vmargin = 1 in]{geometry}
\usepackage{hyperref}
\usepackage{tikz-cd}
\usepackage[all]{xy}
\usepackage{enumitem}

\newenvironment{sbm} %Small bracketed matrix -- good for 3x3 and smaller matrices occurring inside complicated formulas. 
    {\left[ \begin{smallmatrix}
    }
    { 
     \end{smallmatrix} \right]
    }
\newcommand{\defn}[1]{\emph{#1}}
\newcommand{\newword}[1]{\emph{#1}}
\newcommand{\into}{\hookrightarrow}

\newcommand{\Gysin}{\mathrm{Gysin}}
\newcommand{\Filt}{\mathrm{Filtration}}

\DeclareMathOperator{\Hom}{Hom}

\DeclareMathOperator{\Spec}{Spec}

\newcommand{\Alt}{\bigwedge\nolimits}

\renewcommand{\AA}{\mathbb{A}}

\newcommand{\CC}{\mathbb{C}}

\newcommand{\FF}{\mathbb{F}}
\newcommand{\GG}{\mathbb{G}}

\newcommand{\QQ}{\mathbb{Q}}

\newcommand{\ZZ}{\mathbb{Z}}

\newcommand{\cA}{\mathcal{A}}
\newcommand{\cB}{\mathcal{B}}

\newcommand{\cD}{\mathcal{D}}

\newcommand{\cF}{\mathcal{F}}

\newcommand{\cI}{\mathcal{I}}

\newcommand{\cO}{\mathcal{O}}

\newcommand{\cY}{\mathcal{Y}}

\newtheorem{theorem}{Theorem}

\newtheorem{thm}[theorem]{Theorem}
\newtheorem{Theorem}[theorem]{Theorem}
\newtheorem{proposition}[theorem]{Proposition}
\newtheorem{prop}[theorem]{Proposition}
\newtheorem{corollary}[theorem]{Corollary}
\newtheorem{cor}[theorem]{Corollary}

\newtheorem{lemma}[theorem]{Lemma}
\newtheorem{lem}[theorem]{Lemma}

\newtheorem{remark}[theorem]{Remark}

\numberwithin{theorem}{section}
\newtheorem{example}[theorem]{Example}

\newcommand{\GSV}{\mathrm{GSV}}
\DeclareMathOperator{\dlog}{\mathrm{dlog}}

\DeclareMathOperator{\gr}{\mathrm{gr}}

\def\tB{\tilde B}

\def\x{\mathbf{x}}
\def\y{\mathbf{y}}

\def\Id{\mathrm{Id}}

\def\oG{\vec{\Gamma}}

\def\Aut{{\rm Aut}}
\def\tG{{\tilde G}}
\def\frozen{\mathrm{frozen}}
\def\tH{\widetilde{H}}
\def\bI{{\bar I}}
\def\prin{{\rm prin}}
\def\tG{{\tilde G}}

\title{Cohomology of cluster varieties. II. \\ Ayclic case}
\author[Thomas Lam]{Thomas Lam}
\address{Department of Mathematics, 
University of Michigan, Ann Arbor, MI 48109, USA.}
\email{tfylam@umich.edu}
\thanks{T.L. was supported by NSF grants DMS-1160726, DMS-1464693, and DMS-1953852. D.E.S. was supported by DMS-1600223, DMS-1855135 and DMS-1854225.}
\author[David E Speyer]{David E Speyer}
\email{speyer@umich.edu}
\begin{document}
\begin{abstract} 
In previous work, we initiated the study of the cohomology of locally acyclic cluster varieties.  In the present work, we show that the mixed Hodge structure and point counts of acyclic cluster varieties are essentially determined by the combinatorics of the independent sets of the quiver.  We use this to show that the mixed Hodge numbers of acyclic cluster varieties of really full rank satisfy a strong vanishing condition.
\end{abstract}
\maketitle
\setcounter{tocdepth}{1}
\tableofcontents

\section{Introduction}
Cluster algebras (skew-symmetric of geometric type) are determined by an integer $(n{+}m) \times n$ extended exchange matrix $\tB$ whose top $n \times n$ part $B$ is a skew-symmetric integer matrix.  The matrix $B$ determines a directed graph $\oG$ on $[n]:=\{1,2,\ldots,n\}$ with edges $i \to j$ whenever $B_{ij} > 0$.  A cluster variety $\cA(\tB)$ is acyclic if for some cluster, the quiver has no directed cycles.  Locally acyclic cluster varieties, introduced by Muller~\cite{Mul}, have an open covering by acyclic cluster varieties.  

In a previous paper~\cite{LS}, we studied the cohomology, the mixed Hodge structure, and point counts of locally acyclic cluster varieties.  In particular, we proved a curious Lefschetz theorem for Louise cluster varieties of full rank.  In this paper, we specialize to the case of acyclic cluster varieties and obtain stronger and more explicit results.

\subsection{Main results}
Let $\cA$ be an acyclic cluster variety of full rank with $n$ mutable variables and $m$ frozen variables.
We put $d = n + m$, the dimension of $\cA$.  The cohomology $H^{\ast}(\cA):=H^*(\cA,\CC)$ is of mixed-Tate type, meaning that the mixed Hodge numbers $\dim H^{k, (s,t)}(\cA)$ are zero for $s \neq t$ (see Section~\ref{sec:MH} and \cite{LS}).
So $H^k(\cA)$ has a Deligne splitting\footnote{We work with cohomology with complex coefficients throughout, but we mention \cite[Theorem 8.3]{LS}: the Deligne splitting is defined over $\QQ$, in the sense that each $H^{k, (s,s)}(\cA)$ has a basis in $H^k(\cA, \QQ)$.} $H^k(\cA, \CC) = \bigoplus_{s } H^{k,(s,s)}(\cA)$. 
The overarching goal of this paper and our research is to obtain explicit descriptions of the spaces $H^{k, (s,s)}(\cA)$.

Our main structural results are (Theorem~\ref{thm:explicitcomplex}) an explicit complex of finite dimensional vector spaces whose cohomology computes $H^{\ast}(\cA)$, complete with its mixed Hodge structure, and (Theorem~\ref{thm:E1}) a filtration on that complex such that the first page of the resulting spectral sequence has an explicit description in terms of the independence complex of the initial quiver. 
Our main applications are a cohomology vanishing result (Theorem~\ref{thm:main}) and concrete computations of low degree cohomology groups (Section~\ref{sec:small}).
 
In the remainder of the introduction, we assume for simplicity of exposition that $\cA$ is of \textbf{\newword{really full rank}}. This condition means that the rows of the extended exchange matrix $\tB$ span $\ZZ^n$ over $\ZZ$.  In the body of the paper, we will discuss the more general full rank case, where we only impose that the rows span $\QQ^n$ over $\QQ$.  The reduction of the full rank case to the really full rank case is discussed in Section~\ref{sec:decomposition}.

\subsection{The anticlique stratification}
We write $\vec{\Gamma}$ for the directed quiver of the initial $B$-matrix: $\vec{\Gamma}$ has $n$ vertices and has an edge from vertex $i$ to vertex $j$ if $B_{ij}>0$.
We write $\Gamma$ for the undirected graph underlying $\vec{\Gamma}$. 
Let $x_1$, \dots, $x_n$, $x_{n+1}$, \dots, $x_{n+m}$ be the cluster variables in the initial seed, with $x_1$, \dots, $x_n$ mutable and $x_{n+1}$, \dots, $x_{n+m}$ mutable.

An \defn{anticlique} or \defn{independent set} of a graph $G$ is a subset $I$ of the vertices of $G$ such that there are no edges between vertices of $G$.
Let $\cI$ denote the set of anticliques of $\Gamma$ and let $\cI_k \subset \cI$ denote those anticliques of size $k$.  
For an anticlique $I$ of $\Gamma$, let $\cO_I$ be the locally closed subvariety of $\cA$ where the initial cluster variable $x_i$ is zero for $i \in I$ and nonzero for $i \not\in I$.
In Section~\ref{sec:stratification}, we show that $\cA = \bigsqcup_I \cO_I$.  With our really full rank assumption, for $I$ an anticlique of size $k$, each $\cO_I$ is the product of an affine space of dimension $k$ and a torus of dimension $n+m-2k$ (Proposition~\ref{prop:cO}).
%\footnote{This statement and those in the rest of this section has been simplified by our introductory assumption that $\cA$ has really full rank; see Section~\ref{sec:stratification} for the full rank case.}
We obtain the elegant formula (Proposition~\ref{prop:pointcount})
\[
\# \cA(\FF_q) = \sum_k |\cI_k| q^k (q-1)^{n+m-2k}.
\]

The cohomology $H^*(\cO_I)$ is an exterior algebra on a vector space of dimension $n+m-2k$; this vector space is spanned by the classes in $H^1$ whose de Rham representatives are $\dlog x_i := \tfrac{d x_i}{x_i}$ for $i \not\in I$.  We define $\alpha_j := \sum_{i=1}^{n+m} \tB_{ij} \dlog x_i$.
Let $G^0$ be the exterior algebra on the basis $\dlog x_i$ and, for an anticlique $I$, let $G^I$ be the $\CC \langle \dlog x_i : i \not\in I \rangle$-submodule of $G^0$ generated by $\bigwedge_{i \in I} \alpha_i$.
We have (Proposition~ \ref{prop:HOI}) a natural isomorphism $H^{\ast}(\cO_I) \cong (G^I)^{\ast +k}$, where $k = |I|$.

Let $J$ be an anticlique, $j$ an element of $J$ and $I = J \setminus \{ j \}$. 
We define a map $\rho^J_I : G^I \to G^J$ as follows.  For $\theta \in G^I$, write $\theta = \theta_1 + \theta_2 \wedge \dlog{x_j}$ where $\theta_1$ and $\theta_2$ do not involve $\dlog x_j$.
We set $\rho^J_I(\theta) = \theta_2 \wedge \alpha_j$.
We will give an explicit basis for $G^I$ in Section~\ref{sec:basis}, and explain how to write the maps $\rho^J_I$ in that basis.

Set $G^k = \bigoplus_{I \in \cI_k} G^I$. The direct sums of the maps $\rho^J_I$ define a complex $G^0 \to G^1 \to G^2 \to \cdots$, which we denote $G^{\bullet}$. The $\rho^J_I$ maps are graded of degree $0$, so $G^{\bullet}$ splits as a direct sum $\bigoplus_s G^{\bullet,s}$ of subcomplexes $G^{\bullet,s}$.  We call $G^\bullet$ the \emph{Gysin complex}.
The following result gives an explicit description of the cohomology and mixed Hodge structure of $\cA$.
\begin{Theorem}[=Theorem~\ref{thm:explicitcomplex}]
 \label{thm:mainexplicitcomplex}Let $\cA$ be an acyclic cluster algebra of really full rank. With the notation above, we have $H^{k, (s,s)}(\cA) \cong H^{k-s}(G^{\bullet, s})$. (We also recall that $H^{k, (s,t)}(\cA)=0$ for $s \neq t$).
\end{Theorem}

\subsection{The anticlique spectral sequence for principal coefficients}
The \defn{principal coefficients} extended exchange matrix $\tB =B_\prin$ is the extended exchange matrix with $m=n$ such the bottom half of $B_\prin$ is the $n \times n$ identity matrix.
As we discuss in Section~\ref{PrincipalReduce}, the really full rank case can be reduced to the principal coefficients case.

For a graph $G$, the \defn{anticlique (or independence) complex} of $G$ is the simplicial complex $\cI(G)$ whose vertices are the vertices of $G$ and whose faces are the anticliques.
For a subset $S$ of the vertices of $\Gamma$, we write $\Gamma_S$ for the induced subgraph of $\Gamma$ on the vertices of $S$.

\begin{theorem}[=Theorem~\ref{thm:E1}] \label{mainE1} Let $\cA = \cA(B_\prin)$ be an acyclic cluster variety with principal coefficients.  There is a descending filtration  $F^\bullet G^\bullet$ of the complex $G^\bullet$ such that the associated graded complex decomposes as a direct sum over subcomplexes indexed by arbitrary subsets $D$, $E \subseteq [n]$:
\[ \gr G^\bullet = \bigoplus_{(D,E)} \gr G^\bullet(D,E). \]
 We have $H^*(\gr G^\bullet(D,E)) \cong \tH^{*-1}(\cI(\Gamma_{E {\setminus} D}))$, the reduced cohomology of an independence complex.
\end{theorem}

We then use properties of the cohomologies $\tH^{*-1}(\cI(\Gamma_{E {\setminus} D}))$ to study the cohomology of $\cA$.  The filtration $F^\bullet G^\bullet$ gives rise to a spectral sequence $E_{r,\Filt}^{pq}$, and the $E_1$-page has groups given by the cohomologies of the graded pieces of $\gr G^\bullet$.  In Theorem~\ref{thm:E1map}, we describe the differentials $\partial_1$ on the $E_1$-page in terms of natural maps constructed from the independence complexes.

\subsection{Cohomology vanishing}
%\textbf{In this section, we drop the assumption that $\tB$ has principal coefficients, but keep the assumption that $\tB$ has really full rank.}
%We put $d = n+m$, the dimension of $\cA$.

The variety $\cA$ is smooth and affine and of dimension $d = n+m$.  By general results about mixed Hodge structure, this implies that $H^{k,(s,s)}(\cA)$ vanishes unless $0 \leq k \leq d$ and  $\tfrac{1}{2} k \leq s \leq k$.
However, the lower bound for $s$ can be dramatically improved.
\begin{theorem}[=Theorem~\ref{thm:vanishing}] \label{thm:main}
Let $\cA$ be an acyclic cluster variety of really full rank. Then we have $H^{k,(s,s)}(\cA) =0$ unless $0 \leq k \leq d$ and  $\max(\tfrac{2}{3} k, 2k-d) \leq s \leq k$.
%for $3s<2k$ and $d+s < 2k$.
\end{theorem}

In~\cite{LS}, we proved the curious Lefschetz theorem for (even-dimensional) Louise cluster varieties of full rank; this implies that we have the curious Lefschetz symmetry $H^{k,(s,s)}(\cA) \simeq H^{k+d-2s,(d-s,d-s)}(\cA)$ for any acyclic cluster variety of full rank.  This isomorphism swaps the two lower bounds of Theorem~\ref{thm:main}.

The bound of Theorem~\ref{thm:main} does not hold for \emph{locally} acyclic cluster varieties of really full rank. 
Consider the $20$-dimensional cluster variety with the same cluster type as the Grassmannian ${\rm Gr}(5,11)$ and no frozen variables.
By~\cite{GL}, the mixed Hodge structure of this cluster variety is encoded by the $(q,t)$-Catalan number $C_5(q,t) = (q^{10}+ q^9 t + \cdots + q t^9 + t^{10}) + (q^8 t + q^7 t^2 + \cdots + q^2 t^7 + q t^8) + (q^7 t + 2 q^6 t^2 + 2 q^5 t^3 + 2 q^4 t^3 + 2 q^3 t^5 + 2 q^2 t^6 + q t^7) + (q^6 t + q^5 t^2 + 2 q^4 t^3 + 2 q^3 t^4 + q^2 t^5 + q t^6) + (q^4 t^2 + q^3 t^3 +q^2 t^4)$. The $q^6 t$ term reflects a nonzero contribution in $H^{8, (5,5)}$, which violates the bound of Theorem~\ref{thm:main}. 

More generally, the monomial $q^{\binom{n-1}{2}} t$ occurs in $C_n(q,t)$. We can see this using Haglund's formula~\cite{Hag} for $C_n(q,t)$ in terms of ``bounce" and ``area"; the partition $(n-1)$ has bounce equal to $1$ and area equal to $\binom{n-1}{2}$. This corresponds to having $H^{2n-2, (n,n)} \neq 0$ in a cluster variety of dimension $n(n-1)$.

The following corollary follows from Theorem~\ref{thm:main} and curious Lefschetz symmetry.
\begin{cor} \label{cor:vol}
We have $\dim H^d(\cA) = 1$ and $H^d(\cA) \simeq H^{d,(d,d)}(\cA)$. 
\end{cor}
Indeed, it follows from the results of~\cite{LS} that for any cluster torus $T \subseteq \cA$, the map $H^{d, (d,d)}(\cA) \to H^d(T)$ is an isomorphism.
In fact, if $\{ x_1, x_2, \ldots, x_d \}$ is the cluster corresponding to a cluster torus $T$, then the results of~\cite{LS} show that $\Alt_{i=1}^d \dlog x_i$ extends to a closed differential form on $\cA$, which represents a generator of the one-dimensional space $H^{d,(d,d)}(\cA)$. 
So Theorem~\ref{thm:main} shows that this differential form represents  a generator of $H^d(\cA)$.

%Indeed, for any cluster torus $T \subseteq X$, the map $H^d(X, \CC) \longrightarrow H^d(T, \CC)$ is an isomorphism.
Corollary~\ref{cor:vol} has a number of applications.  Recent work in the theory of scattering amplitudes has involved the study of certain integrals on Grassmannians and on other cluster varieties \cite{book,ABL,AHL}.  On the one hand, Corollary~\ref{cor:vol} says that $\cA$ has a unique (in cohomology) volume form, and this is the one used in integrals computing scattering amplitudes.  On the other hand, it follows from Corollary~\ref{cor:vol} that any top-dimensional contour is homologous (up to torsion) to the natural compact contour $S^d \subset (\CC^\times)^d \cong T \subset \cA$.  Note however that Grassmannians are not acyclic, but are locally acyclic \cite{MS}.

In another direction, Corollary~\ref{cor:vol} is relevant in mirror symmetry, for example to conjectures the cluster varieties are large complex structure limit points.  For the case of open positroid varieties and open Richardson varieties, see \cite[Section 8]{HLZ}, \cite{LT}.

For more applications of the cohomology of cluster varieties, we refer the reader to \cite{LS,GL}.

\subsection{Formulae for some mixed Hodge numbers}
In Section~\ref{sec:small}, we compute from $F^\bullet G^\bullet$ the mixed Hodge groups $H^{k,(s,s)}(\cA)$ for $s \leq 3$.  In particular, we obtain the following results.  Note that in \cite{LS} (see Theorem~\ref{thm:LS}) we have already completely described the top-weight subspace $H^{k,(k,k)}(\cA) \subset H^k(\cA)$.
\begin{theorem} \
\begin{enumerate}
%\item We have $H^{0,(0,0)}(\cA) = H^0(\cA)$ and $H^{1,(1,1)}(\cA) = H^1(\cA)$.
\item The $(2,2)$-part of $H^*(\cA)$ is equal to $H^{2,(2,2)}(\cA) \oplus H^{3,(2,2)}(\cA)$.  We have an isomorphism $\dim H^{3,(2,2)}(\cA) \simeq H^1(\Gamma)$.
\item The $(3,3)$-part of $H^*(\cA)$ is equal to $H^{3,(3,3)}(\cA) \oplus H^{4,(3,3)}(\cA)$.  The group $H^{4,(3,3)}(\cA)$ is the direct sum of $H^{1,(1,1)}(\cA) \wedge H^{3,(2,2)}(\cA)$, and another subspace described in Proposition~\ref{prop:principals=3}.  
%When $\Gamma$ is a tree, we have $\dim(V) = \sum_{i \in [n]} \binom{d_i-1}{2}$, where $d_i$ denotes the degree.
\end{enumerate}
\end{theorem}

In Section~\ref{sec:examples}, we compute the mixed Hodge numbers in the case that $\Gamma$ is a star, and discuss some other examples.

\section{Acyclic cluster varieties}
We use the notation $[n]:=\{1,2,\ldots,n\}$ and $[m,n] := \{ m, m+1, m+2, \ldots, n \}$.
\subsection{Cluster algebras}
An \newword{extended exchange matrix} is a $(n{+}m) \times n$ matrix $\tB = (\tB_{ij})$ such that the top $n \times n$ square submatrix $B$, is skew-symmetric. (We anticipate no difficulty in extending our results to the skew-symmetrizable case, but restrict to the skew-symmetric case for convenience.)  For $k \in \{1,2,\ldots,n\}$, we define the mutation of $\tB$ in the direction $k$ to be the extended exchange matrix $\tB'$, given by
\[ \tB'_{ij} = \begin{cases} - \tB_{ij} & \mbox{if $i=k$ or $j=k$,}\\ 
\tB_{ij} + [\tB_{ik}]_{+} [\tB_{kj}]_{+} - [\tB_{ik}]_{-} [\tB_{kj}]_{-} & \mbox{otherwise,} \\ \end{cases} \]
where $[x]_+ = \max(x,0)$ and $[x]_- = \min(x,0)$.
The \newword{top part} of $\tB$ is the $n \times n$ square matrix $B$ given by $B_{ij} = \tB_{ij}$.

Let $\cF$ be a field isomorphic to $\CC(t_1, \ldots, t_{n+m})$. 
A \newword{seed} $t = (\x, \tB)$ in $\cF$ consists of $n{+}m$ elements $x_1$, $x_2$, \ldots, $x_{n+m}$ generating $\cF$ as a field over $\CC$ and 
an extended exchange matrix $\tB$. 
%The variables $x_j$ for $1 \leq j \leq n$ are called \newword{cluster variables}; the variables $x_{n+1}$, \ldots, $x_{n+m}$ are called \newword{frozen variables}; the $n+m$-tuple $(x_1, \ldots, x_{n+m})$ is called a \newword{cluster}.
Given a seed $(\x,\tB)$ and  an index $k$ between $1$ and $n$, the mutation of $(\x,\tB)$ at $k$ is the new seed $(\x',\tB')$, where
\begin{align}\label{eq:mutation}
x'_k &= \frac{\prod_{i} x_i^{[\tB_{ik}]_+} \ + \ \prod_i x_i^{[-\tB_{ik}]_+}}{x_k} \\ 
x'_i &= x_i  & \mbox{if $i \neq k$.} 
\end{align}
%\[\begin{array}{rcl}
%% \tB'_{ij} &=& \begin{cases} - \tB_{ij} & \mbox{if $i=k$ or $j=k$,}\\ 
%%\tB_{ij} + [\tB_{ik}]_{+} [\tB_{kj}]_{+} - [\tB_{ik}]_{-} [\tB_{kj}]_{-} & \mbox{otherwise,} \\ \end{cases} \\
%x'_j &=& \begin{cases} 
% \frac{\prod_{i} x_i^{[\tB_{ik}]_+} \ + \ \prod_i x_i^{[-\tB_{ik}]_+}}{x_k}  &j=k \\ 
% x_j & \mbox{otherwise.} \\
%\end{cases} \\
%\end{array} \]
%We note that mutation leaves the frozen variables unaltered.

We continue mutating on all possible indices, producing new seeds.
The $(n+m)$-tuples $(x_1, \ldots, x_{n+m})$ produced in this manner are called clusters, and the individual $x_i$ are called cluster variables.  The cluster variables $x_1, x_2, \ldots, x_n$ are called \newword{mutable}.  The cluster variables $x_{n+1}, x_{n+2}, \dots, x_{n+m}$ are the same in every cluster and are called \newword{frozen}; we shall often denote them by $y_1 = x_{n+1}$, $y_2 = x_{n+2}$, \dots, $y_m = x_{n+m}$ as well.  

The $\CC$-subalgebra of $\cF$ generated by all the cluster variables, and the reciprocals of the frozen variables, is the \newword{cluster algebra} $A(\x,\tB)$, or simply $A(\tB)$ or $A$.  
We define the \newword{cluster variety} to be the affine variety $\Spec A$ and denote it by $\cA$ or $\cA(\tB)$ or $\cA(\x, \tB)$.
In general, the cluster algebra need not be finitely generated, but we will describe conditions in Section~\ref{sec:acyclic} under which it is.  We say that cluster algebra $A(\x,\tB)$, or the cluster variety $\cA(\x,\tB)$, has rank $n$.

\subsection{Laurent phenomenon}
The Laurent phenomenon~\cite[Theorem 3.1]{CA1}
states that, for any seed $(\tB, (x_1, \ldots, x_{n+m}))$, the cluster algebra $A$ is contained in the Laurent polynomial ring $\CC[x_1^{\pm}, \ldots, x_{n+m}^{\pm}]$. 
This containment turns into an equality if we invert $x_1 \cdots x_n$, that is, $(x_1 x_2 \cdots x_n)^{-1} A = \CC[x_1^{\pm}, \ldots, x_{n+m}^{\pm}]$. 
Geometrically, this means that the open subset of $\cA$ where $x_1 \cdots x_n \neq 0$ is isomorphic to a $(n+m)$-dimensional torus $(\CC^{\ast})^{n+m}$. 
We call such a torus (for any cluster in $A$) a \newword{cluster torus} in $\cA$. 
We caution that it is almost never true that the cluster variety is the union of the cluster tori. Here is the simplest example:

\begin{example}
Let $\tB = \begin{sbm} 1 \\ 0 \end{sbm}$. Then $A = \CC[x_1, x_2, y^{\pm}] / (x_1 x_2 = y+1)$, with the two clusters $(x_1, y)$ and $(x_2, y)$. So the cluster variety $\cA$ is the hypersurface $\{ x_1 x_2 = y+1 \}$ in $\CC^2 \times \CC^{\ast}$, which is isomorphic to the open locus $x_1 x_2 \neq 1$ in $\CC^2$. The two cluster tori are $\cA \cap \{ x_1 \neq 0 \}$ and $\cA \cap \{ x_2 \neq 0 \}$. We note that the point $(x_1, x_2, y) = (0,0,-1)$ is in neither cluster torus.
\end{example}

By the Laurent phenomenon, the cluster algebra $A$ is contained in $\bigcap \CC[x_1^{\pm}, \ldots, x_{n+m}^{\pm}]$, where the intersection is over all clusters $(x_1, \ldots, x_{n+m})$. 
This intersection is known as the \defn{upper cluster algebra} and denoted $U(\x, \tB)$, $U(\tB)$ or $U$.
We have $A \subseteq U$, but we need not have equality, nor need $U$ be finitely generated~\cite{SpeyerFinGenCounterExample}.

%Finally, we discussed above what happens when we invert $x_1 \cdots x_n$. It is also natural to invert only a subset of these variables.
Let $I$ be a subset of $\{ 1,2,\ldots, n\}$. Let $\tB_{[n] \setminus I}$ be the extended exchange matrix formed by deleting the columns indexed by $I$, and reindexing the rows labeled by $I$ to be frozen.
Finally, let $\x_I$ be the same set of variables $(x_1, \ldots, x_{n+m})$ as in $\x$, but with the mutable variables $\{ x_i \}_{i \in I}$ relabeled as frozen. Then we have
\begin{equation} A(\x_I, \tB_{[n] \setminus I})  \subseteq A(\x,\tB)[x_i^{-1} |i \in I] \subseteq  U(\x,\tB)[x_i^{-1} |i \in I] \subseteq U(\x_I,\tB_{[n] \setminus I}) . 
\label{AUChain}
\end{equation}
Geometrically, $\Spec  A(\x,\tB)[x_i^{-1} |i \in I]$ is the open subset $\{ x \in \cA \mid x_i \neq 0 \;{\rm for} \; {\rm all }\; i \}$ in $\cA$; we denote this open subset as $\cA_I$. 

\subsection{Acyclicity} \label{sec:acyclic}

Let $t = (\x,\tB)$ be a seed.  Let $\oG(t) = \oG(\tB)$ be the directed graph with vertices $\{1,2,\ldots,n\}$ and a directed edge $i \to j$ whenever $\tB_{ij} > 0$.   We write $\Gamma(t)$ for the underlying undirected graph of $\oG(t)$.  We say that a seed $(\x,\tB)$ is \defn{acyclic} if $\oG(t)$ has no oriented cycles; see \cite{CA3}.  We say that $\cA$ is acyclic if some seed of $\cA$ is acyclic.  We record a number of results of Berenstein, Fomin and Zelevinsky~\cite{CA3} and Muller~\cite{Mul} concerning acyclic cluster algebras.  

\begin{theorem}[{\cite[Theorem 1.18]{CA3}}]
Let $A(\tB)$ be an acyclic cluster algebra.  Then the cluster algebra coincides with the upper cluster algebra: $A(\tB) = U(\tB)$.
\end{theorem}

The following result states that the open subsets $\cA_I \subset \cA$ of an acyclic cluster variety are themselves cluster varieties.
\begin{theorem}[{\cite[Lemma 3.4]{Mul}}] \label{thm:Mul}
Let $(\x,\tB)$ be an acyclic seed of a cluster algebra $A(\tB)$, and suppose $I \subseteq \{1,2,\ldots,n\}$.  
Then all the containments of~\eqref{AUChain} are equalities.
%Then we have $ A(\x_I,\tB^I)  =A(\x,\tB)[x_i^{-1} |i \in I]$.
\end{theorem}

\begin{theorem}[{\cite[Corollary 1.21]{CA3}}] \label{BigCA3Theorem}
Let $(\x, \tB)$ be an acyclic seed of a cluster algebra $A(\tB)$.  Index the elements of $\x$ as $(x_1, \ldots, x_{n+m})$, with $x_{n+1}, x_{n+2}, \ldots, x_{n+m}$ frozen. For $1 \leq j \leq n$, let $x'_j$ be the variable obtained by mutating at $x_j$. Then $\cA$ is isomorphic to the subvariety of $\CC^{2n} \times (\CC^{\ast})^m$ cut out by the equations
 \begin{equation}\label{eq:exchange}
  x_j x'_j = \prod_{i=1}^{m+n} x_i^{[\tB_{ij}]_+} +  \prod_{i=1}^{m+n} x_i^{[-\tB_{ij}]_+} . %\label{BigCA3Eqns}
  \end{equation}
 Here $x_j$ and $x'_j$, for $1 \leq j \leq n$ are coordinates on $\CC^{2n}$ and $x_j$ for $n+1 \leq j \leq n+m$ are coordinates on $(\CC^{\ast})^m$.
\end{theorem}

We will say that an edge $a \to b$ of $\oG(t)$ is a \newword{separating edge} if there does not exist any bi-infinite directed path \ldots, $i_{-2}$, $i_{-1}$, $i_0$, $i_1$, $i_2$, \ldots in $\oG(t)$ with $i_0=a$, $i_1=b$. In an acyclic graph, every edge is separating.
Muller noticed the following crucial lemma.
\begin{lemma}[{\cite[Corollary 5.4]{Mul}}] \label{lem:Muller}%Add citation
If $a \to b$ is a separating edge, then there is no point of $\cA$ where $x_a = x_b = 0$.  Thus $\cA = \cA_{\{a\}} \cup \cA_{\{b\}}$. %there is no point of $\cA$ where $x_a = x_b = 0$.
\end{lemma}
%In other words, we can write $\cA$ as $\cA_{\{a\}} \cup \cA_{\{b\}}$. 

Associated to each separating edge $a \to b$, there is a class $\epsilon_{ab}$ in $H^{3, (2,2)}$ which we will call the \newword{edge class}, defined as follows.
Recall that $\cA_{\{a\}}$ and $\cA_{\{ b \}}$ are the open sets where the cluster variables $x_a$ and $x_b$, respectively, are nonzero, and $\cA_{\{ a,b \}} = \cA_{\{a\}} \cap \cA_{\{b\}}$. 
Lemma~\ref{lem:Muller} shows that $X = \cA_a \cup \cA_b$, so we have a Mayer-Vietores sequence $\dots \to H^k(X) \to H^k(\cA_a) \oplus H^k(\cA_b) \to H^k(\cA_{\{a,b\}}) \overset{\delta}{\longrightarrow} H^{k+1}(\cA) \to \dots$. 
Let $\epsilon_{ab} = \delta(\dlog x_a \wedge \dlog x_b)$. 
The boundary map $\delta$ preserves mixed Hodge degree, and $\dlog x_a \wedge \dlog x_b$ is in mixed Hodge degree $(2,2)$ by~\cite[Lemma 2.6]{LS}, so $\epsilon_{ab} \in H^{3, (2,2)}(\cA)$. 

The edge classes provide our first examples of classes lying in $H^{k, (s,s)}$ for $s<k$, and we will eventually see that they are, in a sense, the extremal such classes.
Of course, it is not clear whether $\epsilon_{ab}$ is $0$ but we will eventually see that, if $\oG$ is acyclic, then the edge classes span a subspace of $H^{3, (2,2)}(\cA)$ isomorphic to $H^1(\Gamma)$. 

%
%
%
%Our main theme in this paper is how to use the decomposition of $\cA$ into $\cO_I$ to study $\cA$. 
%First, however, we describe other results of Berenstein, Fomin and Zelevinsky~\cite{CA3} and Muller~\cite{LocAcyc} concerning acyclic cluster algebras.
%
%
%%Address the relative primality hypothesis from BFZ
%%Add citations
%
%
%
%Geometrically, this theorem tells us that the open subset of $\cA$ where $\prod_{i \not\in I} x_i$ is nonzero is itself a cluster variety.

\subsection{Cluster varieties of full and really full rank}

We say that the cluster variety $\cA(\tB)$ or the cluster algebra $A(\tB)$ has \defn{full rank} if the matrix $\tB$ has full rank (that is, rank $n$).  The rank of the matrix $\tB$ is mutation-invariant.  An important special case of full rank exchange matrices is the following.  We say that $\tB$ has \defn{principal coefficients} if $\tB = B_\prin :=\left( \begin{smallmatrix} B \\ \mathrm{Id}_n \end{smallmatrix} \right)$.
% and $d$-principal coefficients if $\tB =  \left( \begin{smallmatrix} B \\ d\,\mathrm{Id}_n \end{smallmatrix} \right)$.

It has long been recognized in the subject of cluster algebras that the rank of the exchange matrix $\tB$ is an important invariant, and the situation where $\tB$ has full rank is particularly important. 

\begin{proposition}[\cite{Mul}]
Suppose that $\tB$ is acyclic and of full rank.  Then the cluster variety $\cA$ is smooth.
\end{proposition}

In~\cite{LS}, we argued that it was particularly natural to consider the stronger condition that the rows of $\tB$ span $\ZZ^n$ as a $\ZZ$-module.  

\begin{prop} %add citation
The cokernel $\ZZ^n / \tB^T \ZZ^{n+m}$, as an abstract finite group, is a mutation invariant of $\tB$. In particular, the property that this cokernel is $0$ is mutation invariant.
\end{prop}

We define a cluster algebra $A$ to be \newword{really full rank} if the rows of each $\tB$ matrix span $\ZZ^n$ over $\ZZ$.
We note that this is equivalent to say that $\ZZ^{n+m}/\tB \ZZ^n$ is torsion free.  This property is clearly preserved by deleting columns of $\tB$, that is to say, by declaring additional variables to be frozen.

\section{Stratification of an acyclic cluster variety}\label{sec:stratification}

\subsection{Stratification}
Fix an $(n+m) \times n$ acyclic extended exchange matrix $\tB$ of full rank with undirected quiver $\Gamma = \Gamma(\tB)$.  Let $A = A(\tB)$ denote the cluster algebra, and $\cA = \cA(\tB)$ the cluster variety.  The cluster variables are denoted $x_1,x_2,\ldots,x_n$, and the frozen variables $y_1 = x_{n+1},\ldots,y_m = x_{n+m}$.  The mutation of $x_i$ is denoted $x'_i$.

Recall that a subset $S \subset V(G)$ of vertices of a graph $G$ is an \defn{independent set} or an \defn{anticlique} if the induced subgraph on $S$ has no edges.  Let $\cI_k$ denote the set of anticliques of size $k$ in $\Gamma$ and $\cI$ denote the set (indeed, simplicial complex) of all anticliques.  Recall that if $\oG(t)$ is acyclic then every edge is a separating edge.  The following result follows immediately from Lemma~\ref{lem:Muller}.

\begin{lemma} \label{lem:anticlique}
Let $t = (\x, \tB)$ be an acyclic seed. For any point $x = (x_1,\ldots,x_{n+m})$ of $\cA$, the set of indices $I \subseteq \{ 1,\ldots, n \}$ for which $x_i = 0$ forms an anticlique in $\Gamma(t)$.
\end{lemma}

For $I \in \cI$, we define 
$\cO_I \subset \cA$ to be the relatively open set where $x_i = 0$ for $i \in I$ and $x_i \neq 0$ for $i \not \in I$.

\begin{corollary}\label{cor:strata}
For an acyclic seed $(\x,\tB)$, we have $\cA = \bigsqcup_{I \in \cI} \cO_I$. 
\end{corollary}

Define $\cD_I := \bigcup_{J \supseteq I,\ J \in \cI} \cO_J$ to be the closure of $\cO_I$ in $\cA$.  The correspondence $I \to \cD_I$ is inclusion-reversing.
  The open set $\cA_I$ introduced earlier is given by $\cA_I = \bigcup_{J \subseteq \bI,\ J \in \cI} \cO_J$, where $\bI = [n] {\setminus} I$.
  % The correspondence $I \to \cA_I$ is inclusion-preserving, and the correspondence $I \to \cD_I$ is inclusion-reversing.

\subsection{Automorphisms of the cluster variety}
Let $\tB$ be an exchange matrix which is full rank but not necessarily really full rank.
So $(\tB \QQ^n \cap \ZZ^{m+n}) / \tB \ZZ^n$ is a finite abelian group, which we denote $X^{\ast}$. %Would we like to jut switch to $X$?
As we will explain in this section, the cohomology $H^{\ast}(\cA)$ has a natural splitting indexed by the elements of $X^{\ast}$.
If $\tB$ has really full rank, then the group $X^{\ast}$ is trivial and all the results in this section are trivial as well.

We define a \newword{cluster automorphism} of $A$ to be a $\CC$-algebra automorphism $\phi$ such that, for each cluster variable $x$, there is a scalar $\zeta(x) \in \CC^{\ast}$ such that $\phi(x) = \zeta(x) x$; we denote the group of cluster automorphisms by $\Aut(A)$.
It follows from the Laurent phenomenon that a cluster automorphism is determined by its values on any cluster $(x_1, \ldots, x_{n+m})$. 
\begin{proposition}[{\cite[Proposition~5.1]{LS}}] \label{prop:AutA} %Check reference
We have an isomorphism  
\[ \Aut(A) \cong \mathrm{Hom}(\ZZ^{n+m}/\tB \ZZ^n, \CC^{\ast}) . \]
More concretely, for $(\zeta_1, \zeta_2, \ldots, \zeta_{n+m}) \in (\CC^{\ast})^{n+m}$, there is a cluster automorphism with $\phi(x_i) = \zeta_i x_i$ if and only if the homomorphism $\ZZ^{n+m} \to \CC^{\ast}$ sending the $i$-th basis vector to $\zeta_i$ factors through $\ZZ^{n+m}/\tB \ZZ^n$.
\end{proposition}

The cohomology $H^{\ast}(\cA, \CC)$ therefore decomposes into eigenspaces for the characters of $\Aut(A)$. 
Moreover, the connected component of the identity in $\Aut(A)$ must act trivially on $H^{\ast}(\cA)$, so only locally connected characters occur.
From Proposition~\ref{prop:AutA}, we deduce:
\begin{cor}
The group of locally constant characters of $\Aut(A)$ is isomorphic to
\[ X^*:= (\tB\QQ^n \cap \ZZ^{n+m})/\tB \ZZ^n . \]
Specifically the pairing between $\Aut(A)$ and $X^{\ast}$ is
\[ \langle (\zeta_1, \zeta_2, \ldots, \zeta_{n+m}),\ (k_1, k_2, \ldots, k_{n+m} ) \rangle = \prod \zeta_j^{k_j} . \]
\end{cor}

For $\chi \in X^{\ast}$, we will write $H^{\ast}(\cA)[\chi]$ for the subspace of $H^{\ast}(\cA)$ which transforms by the character $\chi$ of $\Aut(A)$.

Now, let $I$ be an anticlique of $\Gamma$.
Let $\tB_I$ be the submatrix of $\tB$ where we have kept only the columns indexed by $I$, and relabeled the rows indexed by $[n] \setminus I$ as frozen. 
We remark that, since $I$ is an anticlique, the rows of $\tB_I$ indexed by $I$ are $0$.

We set
\[ X^{\ast}(I) := (\tB_I \QQ^{I} \cap \ZZ^{n+m})/\tB_I \ZZ^{I} . \]
The motivation for considering $X^{\ast}(I)$ is the following:
\begin{proposition}\label{prop:cO} For an anticlique $I \in \cI_k$, the subvariety $\cO_I$ is isomorphic to a product of an affine space, a torus, and a finite set $C$ which is a $\Hom(X^*(I), \CC^*)$-torsor:
%comes equipped with a free transitive action of 
\[
\cO_I  \simeq \AA^k \times \Aut(A(\tB_I)) \simeq \AA^k \times \GG_m^{n+m-2k} \times C
\]
\end{proposition}
\begin{proof}
We use the description of $\cA$ from Theorem~\ref{BigCA3Theorem}, and consider $\cO_I$ as a locally closed subset of $\cA$.

In $\cO_I$, we have $x_j \neq 0$ for $j \in [n] {\setminus} I$.  Thus \eqref{eq:exchange} for such $j$ allows us to eliminate the variable $x'_j$.  On the other hand, for $i \in I$, setting $x_i = 0$ the exchange relation \eqref{eq:exchange} becomes 
\begin{equation}\label{eq:kernel}
\prod_{r=1}^{n+m} x_r^{\tB_{ri}} = -1,
\end{equation} and $x'_i$ can take any value.  It follows that $\cO_I$ is isomorphic to the product of $\AA^k$ (with coordinates $\{x'_i, \mid i \in I\}$) and the subvariety $V$ of the torus $\GG_m^{n+m-k}$ (with coordinates $\{x_j \mid j \in [n] {\setminus} I\} \cup \{y_1,\ldots,y_m\}$) cut out by the $k$ equations \eqref{eq:kernel} for $i \in I$.  Since $\tB$ is full rank, so is the restriction $\tB_I$ of $\tB$ to the columns indexed by $I$.  But $I$ is an anticlique so the rows indexed by $I$ in the matrix $\tB_I$ are all 0.  It follows that $\tB_{[n+m]\setminus I, I}$ has full rank.

We obtain a map of tori $\phi = \phi_{\tB_{[n+m]\setminus I, I}}: \GG_m^{n+m-k} \to \GG_m^k$.  The inverse image of the point $(-1,-1,\ldots,-1) \in \GG_m^k$ is equal to the subvariety $V$.  Since $\tB_{[n+m]\setminus I, I}$ has full rank, the map $\phi$ is surjective, so $V \simeq \ker(\phi)$, the inverse image of the identity $(1,1,\ldots,1)$.  By Proposition~\ref{prop:AutA}, we have $ \ker(\phi) \simeq \Aut(A(\tB_I))$.  The group $\Aut(A(\tB_I))$ is the product of a torus $\GG_m^{n+m-2k}$ (its identity component) with its abelian group $C$ of connected components, $\Hom(X^*(I), \CC^*)$.
%
%The kernel $K = \ker(\phi_{\tB_I})$ is a diagonalizable algebraic group, and is thus the product of a torus $\GG_m^{n+m-2k}$ (its identity component) with its abelian group of connected components.  The group of connected components $K/K_0$ is isomorphic to 
%$$
%K/K_0 \simeq \Hom((\tB_I \QQ^k \cap \ZZ^{[n+m]{\setminus} I})/\tB_I \ZZ^k, \CC^*).
%$$
%The finite abelian group $X^*(I)$ is the group of characters of the finite abelian group $K/K_0$, or equivalently, the group of locally constant chracters on $K$.
\end{proof}

We conclude by proving lemmas about $X^{\ast}(I)$ for future use.

\begin{lemma}
The containment $\tB_I \QQ^{I} \cap \ZZ^{n+m} \subseteq \tB\QQ^{n} \cap \ZZ^{n+m}$ descends to an injection $X^{\ast}(I) \into X^{\ast}$.
\end{lemma}

\begin{proof}
Suppose that $(a_1, \ldots, a_{n+m})$ is an element of $\tB_I \QQ^{k} \cap \ZZ^{n+m}$ which maps to $0$ in the quotient $X^{\ast}$.
Then $(a_1, \ldots, a_{n+m})$ is a $\QQ$-linear combination of the columns of $\tB_I$, and is also a $\ZZ$-linear combination of the columns of $\tB$.
But $\tB$ has full rank, so the column of $\tB$ are linearly independent; we deduce that  $(a_1, \ldots, a_{n+m})$ is a $\ZZ$-linear combination of the columns of $\tB_I$, which means that $(a_1, \ldots, a_{n+m})$ is already $0$ in $X^{\ast}(I)$.
\end{proof}

Thus, from now on, we will consider $X^{\ast}(I)$ as a subgroup of $X^{\ast}$. 

\begin{cor} \label{cor:characters}
Let $I$ and $J$ be two anticliques. We have $X^{\ast}(I) \subseteq X^{\ast}(J)$.
\end{cor}

\begin{proof}
We clearly have $\tB_I \QQ^{I} \cap \ZZ^{n+m} \subseteq \tB_J \QQ^{J} \cap \ZZ^{n+m} \subseteq \tB \QQ^n \cap \ZZ^{n+m}$, inducing maps $X^{\ast}(I) \to X^{\ast}(J) \to X^{\ast}$. 
Since the maps $X^{\ast}(I) \to X^{\ast}$ and  $X^{\ast}(J) \to X^{\ast}$ are inclusions, the map $X^{\ast}(I) \to X^{\ast}(J)$ is also an inclusion and, identifying $X^{\ast}(I)$ and $X^{\ast}(J)$ with subgroups of $X^{\ast}$, we have $X^{\ast}(I) \subseteq X^{\ast}(J)$.
\end{proof}

When necessary, we will denote the inclusion $X^{\ast}(I) \subseteq X^{\ast}(J)$ as $\rho^J_I$.

\begin{lemma} \label{lem:MinimalJg}
Let $\chi \in X^{\ast}$ and let $z = (z_1,\ldots,z_{n+m}) \in \tB \QQ^n \cap \ZZ^{n+m}$ be a lift of $\chi$. 
Write $z = \tB u$ for $u \in \QQ^n$ and set $J(\chi) := \{ i  \mid u_i \not\in \ZZ \} \subseteq [n]$. Then $J(\chi)$ depends only on $\chi$, and not on the choice of representative $z$.
For an anticlique $I$, we have $\chi \in X^{\ast}(I)$ if and only if $J(\chi) \subseteq I$.
\end{lemma}

\begin{proof}
Since $\tB$ has full rank, the representation of $z$ as $\tB u$ is unique.  If $z'$ is another vector in $\tB \QQ^n \cap \ZZ^{n+m}$ representing $\chi$, then $z' = \tB u'$ and $u \equiv u' \mod \ZZ^n$.  Thus $\{i \mid u_i \notin \ZZ\} = \{i \mid u'_i \notin \ZZ\}$ and we see that $J(\chi)$ depends only on $\chi$, as claimed.

Next, suppose that $\chi \in X^{\ast}(I)$. So $\chi$ has a representative $z$ lying in $\tB_I \QQ^I$. Using this representative, it is obvious that $J(\chi) \subseteq I$.

Finally, suppose that $I$ is an anticlique containing $J(\chi)$.
Choose an arbitrary representative $z$ of $\chi$ and write $z = \tB u$. 
Define $u'$ to be the vector in $\QQ^n$ with $u'_i = u_i$ for $i \in I$ and $u'_i = 0$ for $i \not\in I$.
If $i \not\in I$ then $i \not\in J(\chi)$ so $u_i \in \ZZ$; we thus see that $u \equiv u' \bmod \ZZ^n$. 
So $\tB u'$ is another representative of $\chi$ and we see that $\chi \in X^{\ast}(I)$, as required.
\end{proof}

\subsection{Point counts}
The cluster algebra $\cA$ can be defined over any field.  The stratification of $\cA$ by $\cO_I$ is also defined over any field, and in particular over a finite field $\FF_q$.  From Corollary~\ref{cor:strata} and Proposition~\ref{prop:cO}, we obtain formulae for the point counts of $\cA$.
We first do the case of a variety of really full rank:
%
%\begin{corollary}
%Let $\cA$ be an acyclic cluster variety of full rank $n$ with $m$ frozen variables.  Then 
%$$
%\# \cA(\FF_q) = \sum_k q^k (q-1)^{n+m-2k} \sum_{I \in \cI_k} |X^*(I)|.
%$$
%\end{corollary}
%
%If $\tB$ is really full rank, then $\tB_I$ is really full rank for any $I \in \cI$.  Thus $|X^*(I)|=1$ for all $I \in \cI$ in that case and we obtain:
\begin{proposition}\label{prop:pointcount}
Let $\cA$ be an acyclic cluster variety of really full rank $n$ with $m$ frozen variables.  Then 
\[
\# \cA(\FF_q) = \sum_k |\cI_k| q^k (q-1)^{n+m-2k}.
\]
\end{proposition}

The full rank situation is messier, since we need to determine which of the components of $\cO_I$ are defined over $\FF_q$. 
For simplicity, we will limit ourselves to the case where they are all defined over $\FF_q$.
\begin{proposition} \label{prop:fullrankpointcount}
Let $\cA$ be an acyclic cluster variety of full rank $n$ with $m$ frozen variables.
Let $N$ be the exponent of the groups $X^{\ast}(I)$, in other words, the minimal positive integer such that $g^N=1$ for all $g \in X^{\ast}(I)$ and all $I \in \cI$. Suppose that $q \equiv 1 \bmod 2N$. Then
\[
\# \cA(\FF_q) = \sum_k q^k (q-1)^{n+m-2k} \sum_{I \in \cI_k} |X^*(I)|.
\]
\end{proposition}

\begin{proof}[Proof of Propositions~\ref{prop:pointcount} and~\ref{prop:fullrankpointcount}]
Let $I$ be an anticlique with $k$ elements. 
We will show that, under the hypothesis of the Propositions, the stratum $\cO_I$ has $|X^*(I)| q^k (q-1)^{n+m-2k}$ points defined over $\FF_q$; summing on $I$ then proves the result.
Proposition~\ref{prop:cO} says that $\cO_I$ is $C \times \AA^k \times \GG_m^{n+m-2k}$, where $C$ is a principal homogenous space for $X^{\ast}(I)$. The second and third factors have $q^k$ and $(q-1)^{n+m-2k}$ points respectively, so what remains is to show that all the points of $C$ are defined over $\FF_q$. Looking at the proof of Proposition~\ref{prop:cO}, these points are computed by finding $N$-th roots of $-1$, so the hypothesis that $q \equiv 1 \bmod N$ means that they are all defined over $\FF_q$. In the case that $X^{\ast}(I)$ is trivial, $C$ is just a single point and is defined over $\FF_q$.
\end{proof}

%\textbf{Thomas: I started working on it, but it seems pretty messy. I suspect the answer is the formula commented out below, but I didn't finish a proof. 
%Is it worth doing?}

%\begin{lemma}
%Let $q$ be a prime power. The number of solutions to $x^d=-1$ in $\FF_q$ is $\GCD(q-1,2d) - \GCD(q-1,d)$. 
%\end{lemma}
%
%\begin{proof}
%Since $\FF_q^{\times}$ is cyclic, the number of solutions to $x^c=1$ in $\FF_q$ is $\GCD(q-1,c)$. We can count solutions to $x^d=-1$ by counting solutions to $x^{2d}=1$ and subtracting off the solutions to $x^d=1$.
%\end{proof}
%
%Let $X^{\ast}(I) \cong \prod_i \ZZ/d_i(I) \ZZ$ where $d_1(I)$ divides $d_2(I)$ divides $d_3(I)$ and so forth. With this notation, we show:
%\begin{proposition}\label{prop:pointcountfullrank}
%Let $\cA$ be an acyclic cluster variety of full rank $n$ with $m$ frozen variables. Let $q$ be a prime power which is relatively prime to $2 |X^{\ast}|$. Then
%\[ \# \cA(\FF_q) = \sum_k \sum_{I \in \cI_k} q^k (q-1)^{n+m-2k} \prod_i {\big(} \GCD(q-1, 2 d_i(I)) - \GCD(q-1, d_i(I)) {\big)}. \]
%\end{proposition}
%

\subsection{Cohomology of $\cO_I$}

Let $\CC[X^*(I)]$ denote the group ring of $X^*(I)$.  We consider $\CC[X^*(I)]$ a graded ring by placing it completely in degree 0.  

Let 
$$
G^\emptyset:= H^*(\cO_{\emptyset}) = \CC \langle \dlog x_i, \dlog y_i \rangle
$$
be the cohomology of the cluster torus.  It is an exterior algebra over $\CC$ on the $n+m$ generators $\dlog x_i, \dlog y_i$ which have degree $1$.

For $i = 1,2,\ldots,n$, define the one-form $\alpha_i = \sum_r \tB_{ri} \dlog x_r$.
Let $G^I$ be the $\CC\langle \dlog x_i, i \notin I$, $\dlog y_i, i \in [m] \rangle$-submodule of $G^\emptyset$ generated by $\bigwedge_{i \in I}  \alpha_i$.  Equip $G^I$ with the grading inherited from $G^\emptyset$.  Also, let $L^I$ denote the quotient of $L^\emptyset := G^\emptyset$ by the ideal generated by the relations $\alpha_i = 0$ for $i \in I$, and $\dlog x_i = 0$ for $i \in I$.  Since these relations are homogeneous, $L^I$ also inherits a grading from $G^\emptyset$.  Then as modules over $\CC\langle \dlog x_i, i \notin I$, $\dlog y_i, i \in [m] \rangle$, we have an isomorphism $(L^I)^* = (G^I)^{* + k}  $ given by $\theta \mapsto \theta \wedge \bigwedge_{i \in I}  \alpha_i$.

If $J = I \sqcup \{j\}$, there is a graded map $\rho_I^J: L^I \to L^J$ of degree $-1$, defined as follows.  For $\theta \in L^I$, write $\theta = \theta_1 + \theta_2 \wedge \dlog x_j$ where $\theta_1$ and $\theta_2$ do not involve $\dlog x_j$.  Then $\rho_I^J(\theta) := \theta_2$.  Recall from Corollary~\ref{cor:characters}, that we have defined an inclusion $\rho_I^J: X^*(I) \hookrightarrow X^*(J)$.  Let $\rho_I^J: \CC[X^*(I)] \to \CC[X^*(J)]$ also denote the induced map on group algebras.

\begin{proposition} \label{prop:HOI}
For $I \in \cI_k$, we have an isomorphism of graded rings
\begin{equation}\label{eq:isom}
H^*(\cO_I) \simeq L^I \otimes_\CC \CC[X^*(I)]
\end{equation}
where in the first factor the multiplication is wedge product, and in the second factor it is induced by the (commutative) group structure of $X^*(I)$.  

Let $J$ be an anticlique with $k+1$ element, let $j$ be an element of $J$ and write $I = J \setminus \{ j \}$.
Then $\cO_j$ is a hypersurface in $\cO_I \cup \cO_J$ and the residue map  $\rho_I^J: H^*(\cO_I) \to H^*(\cO_J)$ is given by the map
\begin{equation}\label{eq:res}
 \rho_I^J \otimes \rho_I^J: L^I \otimes_\CC \CC[X^*(I)] \to L^J \otimes_\CC \CC[X^*(J)].
\end{equation}
\end{proposition}
%Note that the residue map $\rho_I^J: H^*(\cO_I) \to H^*(\cO_J)$ is defined as the sum of the residue maps $H^*(\cO_I) \to H^*(C)$ over all connected components $C \subset \cO_J$.  But it is convenient to keep track of all the connected components together.

\begin{proof}
According to Proposition~\ref{prop:cO}, we have that $\cO_I$ is isomorphic to $\AA^k \times V$, where $V$ is cut out of the torus   $T = \GG_m^{n+m-k}$ by the equations \eqref{eq:kernel}.  Thus $H^*(\cO_I) \simeq H^*(V)$.  The cohomology $H^*(T)$ is the exterior algebra on generators $\dlog x_i$ for $i \in [n] {\setminus} I$ and $\dlog y_i$ for $i = 1,2,\ldots,m$.  The natural map $H^*(T) \to H^*(V)$ has kernel generated by $\alpha_i$, $i \in I$ (obtained from \eqref{eq:kernel}), and the image can be identified with the cohomology of the identity component $V^0 \subset V$.  It follows that the cohomology $H^*(V^0)$ is isomorphic to $L^I$ as graded rings.

By Proposition~\ref{prop:cO}, the ring $\CC[X^*(I)]$ can be identified with the ring of locally constant functions on $\cO_I$.   The isomorphism \eqref{eq:isom} follows.

For \eqref{eq:res}, we note that the subset $\cO_J$ is cut out of $\cD_I$ by the equation $x_j= 0$.  The description of the residue map $\rho_I^J: L^I \to L^J$ follows immediately.  Finally, note that the identification of connected components of $\cO_I$ with $ \Hom(X^*(I), \CC^*)$ is compatible with the inclusions $\rho_I^J: X^*(I) \to X^*(J)$.
\end{proof}

\begin{remark}
We could define the residue map $\rho_I^{J,a}: H^*(\cO_I) \to H^*(\cO_J^{(a)})$ for each connected component $\cO_J^{(a)} \subset \cO_J$.  For our current purposes, it is convenient to keep track of all the connected components together.
\end{remark}

\section{A basis for $G^I$} \label{sec:basis}

Let $\tB$ be an $(n+m) \times n$ extended exchange matrix of full rank.  For $R \subseteq [n+m]$ and $C \subseteq [n]$, write $\tB_{R,C}$ for the submatrix of $\tB$ with rows $R$ and columns $C$, and write $\tB_C$ for $\tB_{[n+m], C}$.

Let $I$ be a $k$-element anticlique. For $A \subseteq [n+m] \setminus I$, let
\[ \theta(A,I) = \bigwedge_{a \in A} \dlog x_a \wedge \bigwedge_{i \in I} \alpha_i \]
where the wedges are ordered by the induced linear order on $A$ and $I$ as subsets of $[n+m]$.
So $G^I$ is spanned by the $2^{n+m-k}$ wedges $\theta(A,I)$, as $A$ varies over subsets of $[n+m] \setminus I$.

We know, however, that the dimension of $G^I$ is only $2^{n+m-2k}$, so we would like to give a subset of the $\theta(A,I)$ which are a basis for $G^I$. This can be done as follows.  Since $\tB$ has full rank, the columns of $\tB$ indexed by $I$ are linearly independent. Therefore, there must be a set $N(I)$ of $k$ rows of $\tB$ such that $\tB_{N(I), I}$ is invertible. Moreover, since $I$ is an anticlique, $\tB_{II}$ is the $0$-matrix, so $N(I)$ must be disjoint from $I$.
\begin{lemma} \label{lem:Gbasis}
Let $I$ and $N(I)$ be as above. Then the $\theta(A, I)$, for $A \subseteq [n+m] \setminus (I \sqcup N(I))$, is a basis of $G^I$. 
\end{lemma}

\begin{proof}
The condition that $\tB_{I, N(I)}$ is invertible means that the set 
$$ S:= \{ \alpha_i \}_{i \in I} \sqcup \{ \dlog x_j \}_{j \in [n+m] \setminus (I \cup N(I))}$$ is a basis for $\bigoplus_{i \in [n+m] \setminus I} \QQ \dlog x_i$. So wedges of subsets of this set are a basis for the exterior algebra $\Alt^{\bullet} \left( \bigoplus_{i \in [n+m] \setminus I} \QQ \dlog x_i \right)$. We defined $G^I$ as the submodule of this algebra generated by $\bigwedge_{i \in I} \alpha_i$, so a basis for $G^I$ is spanned by wedges of subsets of $S$ containing $\{ \alpha_i \}_{i \in I}$. This is the required result.
\end{proof}

We now compute the maps $\rho^J_I$ in terms of the $\theta(A,I)$. 
\begin{lemma} \label{lem:rhoMapBasic}
Let $J$ be an anticlique and write $J = I \sqcup \{ j \}$. Then
\[ \rho^J_I(\theta(A,I)) = \begin{cases} \pm \theta(A \setminus \{ j \}, J) & j \in A \\ 0 & j \not\in A \end{cases}. \]
\end{lemma}

\begin{proof}
Since $J$ is an anticlique, we have $\tB_{ji}=0$ for all $i \in I$. Therefore, the term $\dlog x_j$ does not occur in $\bigwedge_{i \in I} \alpha_i$. So, if $j \not \in A$, then there is not $\dlog x_j$ in $\theta(A, I)$ and $\rho^J_I(\theta(A,I))=0$. If $j \in A$, then $\theta(A,I) = \pm \dlog x_j \wedge \theta(A \setminus \{ j \}, J)$ and $\rho^J_I(\theta(A,I)) = \pm \alpha_j \wedge \theta(A \setminus \{ j \}, J)  = \pm \theta(A \setminus \{ j \}, J)$, where the $\pm$ signs are independent. 
\end{proof}

Unfortunately, even if $\theta(A,I)$ is in our chosen basis for $G^I$, this does not imply that $\theta(A \setminus \{ j \}, J)$ is. In other words, we can have $A \cap (I \sqcup N(I)) = \emptyset$ but $(A \setminus \{ j \}) \cap N(J) \neq \emptyset$. 
In this case, we need to invert the matrix $\tB_{N(I),I}$ to write $\theta(A \setminus \{ j \}, J)$ as a linear combination of $\theta(A', J)$ for $A'$ disjoint from $J \sqcup N(J)$. 

\begin{remark} \label{rem:WhyPrincipalGood}
The situation is a bit more tractable if there is an injection $\nu : [n] \to [n+m]$ such that we can take $N(I) = \{ \nu(i) : i \in I \}$. In this case, $N(J) = N(I) \sqcup \{ \nu(j) \}$ so the only possible obstacle to having $(A \setminus \{ j \}) \cap N(J)$ empty is if $\nu(j) \in A$. However, in general, such an injection $\nu$ does not exist. Consider the matrix
\[ \tB = \begin{bmatrix}
0&0&0 \\
0&0&0 \\
0&0&0 \\
\hline
1&1&0 \\
1&1&1 \\
0&1&1 \\
\end{bmatrix} .\]
We claim that there is no injection $\nu : [3] \to [6]$ such that $\tB_{\nu(I),I}$ is invertible for every subset $I$ of $[3]$. Using the condition that $\tB_{\nu(i), i} \neq 0$ already reduces us to one of the cases $(\nu(1), \nu(2), \nu(3)) = (4,5,6)$, $(\nu(1), \nu(2), \nu(3)) = (4,6,5)$ or $(\nu(1), \nu(2), \nu(3)) = (5,4,6)$. But, in the second case, $\tB_{\nu(\{ 1,2 \}), \{ 1,2 \}}$ is not invertible; in the third case, $\tB_{\nu(\{ 2,3 \}), \{ 2,3 \}}$ is not invertible and, in the first case, neither $\tB_{\nu(\{ 1,2 \}), \{ 1,2 \}}$ nor $\tB_{\nu(\{ 2,3 \}), \{ 2,3 \}}$ is invertible.
In this case, the troublesome part of $\tB$ occurs in the frozen rows, so we could modify them using Lemma~\ref{lem:MMutablePart}. A more difficult example would be the $6 \times 6$ exchange matrix with this as the first $3$ columns.

There are two natural cases where such an injection $\nu$ exists: The first is the class of a cluster algebra with principal coefficients, defining $\nu(i) = i+n$. We will study this case in detail in Section~\ref{sec:principal}.  The second case is if we have no frozen variables and $\Gamma$ has exactly one perfect matching, such as the $A_{2k}$, $E_6$ and $E_8$ Coxeter diagrams. In this case we can take $\nu : [n] \to [n]$ to be the involution which interchanges the ends of each vertex of the unique matching. We do not study this second case in detail, as there appears to be no way to reduce a general cluster variety to this case, but we remark that such cluster varieties often seem to have unusually small Betti numbers.
\end{remark}

\section{The Gysin complex}\label{sec:cohom}
In this section $\cA$ is an acyclic cluster variety of full rank. 
\subsection{Summary of results from~\cite{LS}} \label{sec:MH}
The cohomology ring $H^*(\cA)$ is equipped with a mixed Hodge structure.  In \cite{LS} we showed that $H^*(\cA)$ is of mixed-Tate (or Hodge-Tate) type, and is thus endowed with a decomposition, called the \newword{Deligne splitting},
$$
H^*(\cA) = \bigoplus_k \bigoplus_{k/2 \leq s \leq k} H^{k,(s,s)}(\cA).
$$
We refer the reader to \cite{LS} for details.  We call $(s,s)$ the \newword{mixed Hodge degree} of the summand $H^{k,(s,s)}(\cA)$, and we say that it has \newword{weight} $2s$.

A \newword{GSV 2-form} \cite{GSV,LS} for an extended exchange matrix $\tB$ is a 2-form
\[
\gamma =  \sum_{i,j} \widehat{B}_{ij} \dlog x_i \wedge \dlog x_j
\]
where $\widehat{B}$ a skew-symmetric $(n+m) \times (n+m)$ matrix whose first $n$ columns equals $\tB$. 
(The initials GSV stand for ``Gekhtman, Shapiro, Vainshtein"; see~\cite{GSV}.)

Let $H^*(\cA)_{st} := \bigoplus_k H^{k,(k,k)}(\cA)$ be the top-weight subring of the cohomology $H^*(\cA)$.  Identify the cohomology $H^*(T)$ of a cluster torus $T \subset \cA$ with the exterior algebra $\CC \langle \dlog x_1,\ldots,\dlog_{x+m} \rangle$.

\begin{theorem}[{\cite{LS}}]\label{thm:LS}
Let $\cA$ be a (locally) acyclic cluster variety of full rank.
\begin{enumerate}
\item The subring $H^*(\cA)_{st}$ can be identified with the subring of $\CC \langle \dlog x_1,\ldots,\dlog_{x+m} \rangle$ consisting of forms that are regular on $\cA$.
\item The subring $H^*(\cA)_{st}$ is generated by the forms $\dlog x_{n+1},\ldots, \dlog x_{n+m}$ and the GSV 2-forms $\gamma_{\Gamma_i} = \gamma_{\tB_{\Gamma_i \cup [n+1,n+m],\Gamma_i}}$ as $\Gamma_i$ varies over the connected components of $\Gamma$.
\item If $\cA$ is connected and has principal coefficients, then the following forms give a basis of $H^*(\cA)_{st}$:
$$
\gamma^j \wedge \bigwedge_{k \in K} \dlog y_k, \qquad j + |K| \leq n
$$
where $\gamma$ is a fixed GSV-form and $K$ varies over the subsets of $[n]$.
\item If $\Gamma$ is a path and $\cA$ has really full rank, then $H^*(\cA) = H^*(\cA)_{st}$.
\end{enumerate}
\end{theorem}

%\subsection{Cohomology}
%Using the isomorphisms $(L^I)^* = (G^I)^{* + k}  $ and $(L^J)^* = (G^J)^{* + k}$, we obtain maps $\rho_I^J: G^I \to G^J$ that are graded of degree $0$, explicitly given as follows.
% For $\theta \in G^I$, write $\theta = \theta_1 + \theta_2 \wedge \dlog x_j$ where $\theta_1$ and $\theta_2$ do not involve $\dlog x_j$.  Then $\rho_I^J(\theta) := \theta_2 \wedge \alpha_j$.  
% 
%Let $G^k:= \bigoplus_{I \in \cI_k} G^I \otimes_\CC \CC[X^*(I)]$, and let 
%$$
%G^\bullet:= G^0 \to G^1 \to \cdots 
%$$
%be the complex whose differential $\partial: G^\bullet \to G^\bullet$ is the direct sum of the $\rho_I^J$.  We write $G^{\bullet,s}$ for the degree $s$ piece of the complex $G^\bullet$.  Since $\partial$ has degree $0$, we have that $G^{\bullet, s}$ is a subcomplex.
%% and $\partial : G^\bullet \to G^\bullet$ for the differential.
%Our main goal for this section is to prove:
%
%\begin{thm}\label{thm:explicitcomplex}
%We have $H^{k,(s,s)}(\cA) \simeq H^{k-s}(G^{\bullet,s})$.
%\end{thm}
%
\subsection{The Gysin spectral sequence}
%Let $Y$ be a smooth algebraic variety and $Z$ a smooth subvariety of codimension one.  Let $U = Y \setminus Z$.  Then there is a Gysin long exact sequence~\cite{PS}:
%\begin{equation}\label{eq:LES}
%\cdots \to H^k(Y) \to H^k(U) \to H^{k-1}(Z) \to H^{k+1}(Y) \to \cdots
%\end{equation}
%where the map $H^k(Y) \to H^k(U)$ is restriction, and $H^k(U) \to H^{k-1}(Z)$ is the residue map.

Let $Y$ be a smooth algebraic variety and $D \subset Y$ a divisor.  We call $D$ a \newword{normal crossings divisor} if each component $D_i$ of $D$ is smooth, and \'etale locally, the intersection $D_{i_1} \cap D_{i_2} \cap \cdots \cap D_{i_k}$ is isomorphic to the intersection of coordinate hyperplanes.  In particular, the intersection $D_{i_1} \cap D_{i_2} \cap \cdots \cap D_{i_k}$ is smooth and of pure codimension $k$.  For the following result see for example \cite[Section 3]{Ara} or \cite{Pet}.

%The following result encodes the repeated application of \eqref{eq:LES} to $Y$ and the $D_i$.

\begin{theorem}%[{\cite{Del}}]
\label{thm:Gysin}
Let $Y$ be a smooth algebraic variety and $D = \bigcup_{i=1}^n D_i \subset Y$ a normal crossings divisor.
We have a spectral sequence $E_r^{pq}(Y,D)$ converging to $H^*(Y)$, and inducing the weight filtration on $H^*(Y)$.  The $E_1$-page is given by
\[
E_1^{pq}(Y,D) =  \bigoplus H^{q-p}{\Big (} (D_{i_1} \cap D_{i_2} \cap \cdots \cap D_{i_p}) {\setminus} \bigcup\nolimits_{i\notin \{i_1,\ldots,i_p\}} D_i {\Big )}
\]
where the sum is over all $p$-tuples $1 \leq i_1 < i_2 < \cdots < i_p \leq n$ of boundary components. The differentials $E_1^{p,q}(Y,D) \to E_1^{p+1,q}(Y,D)$ are the residue maps.
\end{theorem}

We note that, if $D_{i_1} \cap D_{i_2} \cap \cdots \cap D_{i_p}$ is empty, then the corresponding summand is $0$.

If $D = Z$ is a single smooth hypersurface in $Y$, then the spectral sequence of Theorem~\ref{thm:Gysin} reduces to the Gysin long exact sequence of the triple $(Y,U,Z)$:
\begin{equation}\label{eq:LES}
\cdots \to H^k(Y) \to H^k(U) \to H^{k-1}(Z) \to H^{k+1}(Y) \to \cdots
\end{equation}
where $U = Y\setminus Z$.

\subsection{The Gysin complex for an acyclic cluster variety}
\begin{proposition}
Let $\cA$ be a full rank acyclic cluster variety and recall the notation $\cD_{\{ i \}}$ for $\{ x_i =0 \}$. Then $\cD_{\{1 \}}$, $\cD_{\{2 \}}$, \dots, $\cD_{\{n \}}$ is a normal crossings divisor.
%Let $I = \{i_1,i_2,\ldots,i_k\} \in \cI_k$.  Then the divisor $\cD = \cD_{\{ i_1 \}} \cup \cD_{\{ i_2 \}} \cup \cdots \cup \cD_{\{ i_k \}}$ is a normal crossings divisor.
\end{proposition}
\begin{proof}
Being a normal crossings divisor is an \'{e}tale local condition, so let us check it near a particular point $z$ of $\cA$. Let $z$ lie in the stratum $\cO_I$, indexed by an anticlique $I$ and let $k$ be the cardinality of $I$.  Let $\bar I = [n] \setminus I$.
An open neighborhood of $z$ is given by $\cA_{\bar I}$, which is the cluster variety associated to the exchange matrix $\tB_I$.
Since $I$ is an anticlique, the top part of $\tB_I$ is simply the $0$ matrix.
The divisors $\cD_i$ for $i \not\in I$ are disjoint from $\cA_{\bar I}$, so we only need to consider the divisors $\cD_i$ for $i \in I$. 

By \cite[Proposition 5.10]{LS}, we have an \'etale covering map $\eta: \cA' \to \cA_{\bar I}$, where $\cA'$ is a cluster variety with extended exchange matrix $\tB' = \left[ \begin{smallmatrix} 0_{k \times k} \\ d\,\Id_k \\ 0 \end{smallmatrix} \right]$ for some $d \in \ZZ_{>0}$.
The preimage of the divisor $\cD_i$, for $i \in I$, is the hypersurface $x_i = 0$ in $\cA'$. 

But $\cA'$ is simply the product $\cY_d^k \times (\CC^{\ast})^{n+m-2k}$ where $\cY_d = \{ (x,x',y) \mid x x' = y^d+1 \} \subset \CC \times \CC \times \CC^{\ast}$. The functions $x_i$, for $i \in I$, pull back to the $x$-coordinates on the different $\cY_d$ factors. So the equations $x_i=0$ are transverse to each other, and we just need to see that $\{ x= 0 \}$ is a smooth hypersurface in $\cY_d$.
Indeed, in $\cY_d$, the hypersurface $x=0$ is just $d$ copies of the affine line, given by $\{ (0,x', \zeta) : \zeta^d=-1 \}$.
\end{proof}

Thus, we may apply Theorem~\ref{thm:Gysin} to the case of an acyclic cluster variety.
The result is a spectral sequence whose pages we will label $E^{pq}_{r, \Gysin}$.
So
\begin{equation} E^{pq}_{1, \Gysin} = \bigoplus_{1 \leq i_1 < i_2 < \ldots < i_p \leq n} H^{q-p} {\Big (} (\cD_{i_1} \cap \cD_{i_2} \cap \cdots \cap \cD_{i_p}) {\setminus} \bigcup\nolimits_{i\notin \{i_1,\ldots,i_p\}} \cD_i {\Big )} . \label{GysinDumb} 
\end{equation}

We now give a more explicit description of the spaces on the $E_1$ page.
\begin{prop} \label{E1Spaces}
Let $\cA$ be an acylic cluster variety of full rank. In the above notation, we have
\[ E^{pq}_{1, \Gysin} = \bigoplus_{I \in \cI_p} (G^I)^q  \otimes \CC[X^{\ast}(I)] . \]
Here $(G^I)^q$ denotes the $q$-degree part of the graded module $G^I$.
\end{prop}

\begin{proof}
We begin with \eqref{GysinDumb}. 
If $\{ i_1, i_2, \ldots, i_p \}$ is \textbf{not} an anticlique, then $\cD_{i_1} \cap \cD_{i_2} \cap \cdots \cap \cD_{i_p}$ is empty by Lemma~\ref{lem:Muller}, so we may restrict to the summands where $\{ i_1, i_2, \ldots, i_p \}$ is an anticlique.
In that case, $(\cD_{i_1} \cap \cD_{i_2} \cap \cdots \cap \cD_{i_p} ) {\setminus} \bigcup\nolimits_{i\notin \{i_1,\ldots,i_p\}} \cD_i = \cO_I$.
By Proposition~\ref{prop:HOI}, $H^{q-p}(\cO_I)$ is the degree $q-p$ piece of $L^I \otimes \CC[X^{\ast}(I)]$. 
But this is isomorphic to the degree $q$ part of $G^I \otimes \CC[X^{\ast}(I)]$. 
\end{proof}

\begin{remark}
Let $I$ be an anti-clique of size $p$. We note that $G^I$ is supported in degrees $p \leq q \leq m+n-p$. 
We will have $E_{\infty, \Gysin}^{pq} = H^{p+q, (q,q)}(\cA)$.
So, if $H^{k,(s,s)}(\cA)$ is nonzero, then $k = p+q \leq \min(m+n, 2q) = \min(m+n,2s)$ and also $s = q \leq p+q =k$. We thus recover the inequalities $k/2 \leq s \leq k \leq m+n$ which follow directly from the fact that $\cA$ is smooth and affine.\end{remark}

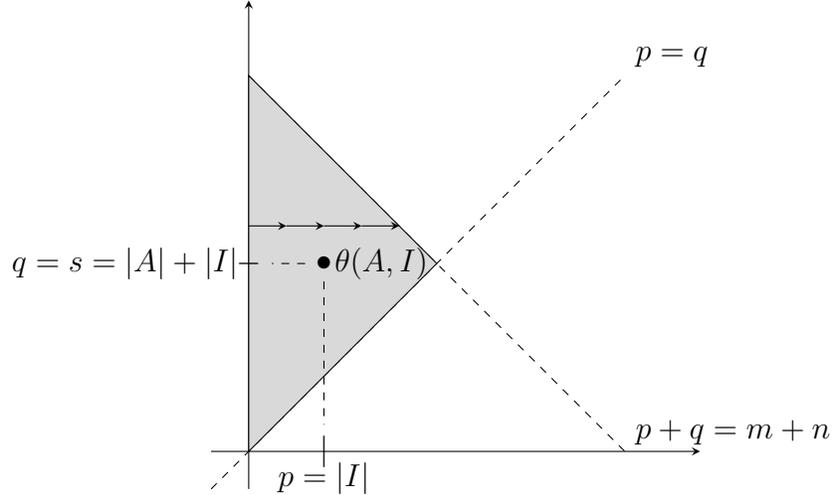
\begin{figure}
\begin{tikzpicture}
\draw[-stealth] (-0.5,0) -- (6,0) ;
\draw[-stealth] (0,-0.5) -- (0,6) ;
\draw[dashed] (-0.5,-0.5) -- (5,5) ;
\draw[dashed] (0,5) -- (5,0) ;
\draw[fill=gray!30] (0,0) -- (0,5) -- (2.5,2.5) -- (0,0);
\node (pTick) at (1,0) {$|$} ;
\node[anchor=north] (pLabel) at (1,0) {$p=|I|$} ;
\node (qTick) at (0,2.5) {$-$} ;
\node[anchor=east] (qLabel) at (0,2.5) {$q=s=|A|+|I|$} ;
\node (point) at (1,2.5) {$\bullet$};
\node[anchor=west] (theta) at (1,2.5) {$\theta(A,I)$};
\draw[dashed] (point) -- (pTick);
\draw[dashed] (point) -- (qTick);
\node[anchor=west] (dimEqn) at (5,0.25) {$p+q=m+n$};
\node[anchor=west] (diagonal) at (5,5.25) {$p=q$};
\draw[-{stealth[scale=2]}] (0,3) -- (0.5,3) ;
\draw[-{stealth[scale=2]}] (0.5,3) -- (1,3) ;
\draw[-{stealth[scale=2]}] (1,3) -- (1.5,3) ;
\draw[-{stealth[scale=2]}] (1.5,3) -- (2,3) ;
\end{tikzpicture}
\caption{The page $E_{1, \Gysin}$. The arrows depict the maps on this page. \label{E1GysinFigure}}
\end{figure}

In Figure~\ref{E1GysinFigure}, the grey shaded triangle indicates where $E_{1, \Gysin}$ is supported.
The horizontal arrows show the direction of the maps on this page.

We next note that we have already computed the maps on the $E_1$ page:
\begin{prop} \label{E1Maps}
With notation as above, the map $E^{pq}_{1, \Gysin} \longrightarrow E^{(p+1)q}_{1, \Gysin}$ is the direct sum of all maps $\rho^J_I  \otimes \rho^J_I : G^I \otimes \CC[X^{\ast}(I)] \to G^J \otimes \CC[X^{\ast}(J)]$, where $J$ is a $p+1$ element anticlique and $I$ is a $p$-element subset of $J$.
\end{prop}

\begin{proof}
The maps on the $E_1$ page of the Gysin spectral sequence are the residue maps $H^{\ast}(\cO_I) \to H^{\ast}(\cO_J)$, and we computed in Proposition~\ref{prop:HOI} that this residue map is given by $\rho^J_I \otimes \rho^J_I$.
\end{proof}

So, for each fixed $q$, we have a complex $E^{0q}_{1, \Gysin} \to E^{1q}_{1, \Gysin} \to \cdots$ and the $E_2$ page of the spectral sequence is the cohomology of this complex.
Our choice of stratification gives an additional simplification.

\begin{prop} \label{E2Collapse}
With notation as above, the Gysin spectral sequence collapses at $E_2$, meaning that all maps on page $E_r$ are $0$ for $r \geq 2$.
\end{prop}

\begin{proof}
The maps in the Gysin spectral sequence are maps of mixed Hodge structures \cite[Section 3]{Ara}, with appropriate shifts: the weight in position $(p,q)$ is shifted by $2p$ (or equivalently, the mixed Hodge degree is shifted by $(p,p)$).   However, $\cO_I$ is the product of a torus, an affine space and a finite set, and $H^k$ of a torus is pure of weight $2k$ (with mixed Hodge degree $(k,k)$). 
So the $(p,q)$ term of the Gysin spectral sequence for the stratification by the $\cO_I$ is of weight $2q$ (and mixed Hodge degree $(q,q)$) in position $(p,q)$. 
So all maps after the first page of the spectral sequence must be $0$.
\end{proof}

We deduce the following description of the cohomology of $\cA$, including its mixed Hodge structure.
In the really full rank case, this is Theorem~\ref{thm:mainexplicitcomplex}.

\begin{theorem} \label{thm:explicitcomplex}
Let
\[ \tG^{ps} = \bigoplus_{I \in \cI_p} (G^I)^s \otimes \CC[X^{\ast}(I)] .\]
where $(G^I)^s$ is the degree $s$ part of $G^I$. Define a map $\tG^{ps} \to \tG^{(p+1)s}$ as the direct sum of the maps $\rho^J_I \otimes \rho^J_I$, where $J$ runs over $p+1$ element anticliques and $I$ runs over $p$-element subsets of $J$. Then the $\tG^{ps}$ form a complex, which we denote $\tG^{\bullet, s}$. We have
\[ H^{k, (s,s)}(\cA) \cong H^{k-s}(\tG^{\bullet, s}) . \]
\end{theorem}

The complex $\tG^\bullet$ has an additional grading by $X^*$ coming from the grading on the factors $\CC[X^*(I)]$.  For $\chi \in X^*$, write $\tG^\bullet[\chi]$ for the summand with degree $\chi$; the subcomplex $\tG^\bullet[\chi]$ is spanned by elements of the form $\theta \otimes \chi$ where $\theta \in G^I$ and $\chi \in X^{\ast}(I)$ for some $I$.  
For the trivial character $\chi_e$, the component $\tG^\bullet[\chi_e]$ is the complex we denoted $G^{\bullet}$ in the introduction; we will  continue to use the notation $G^{\bullet}$ for $\tG^\bullet[\chi_e]$.

When we need to, we will write $\tG^\bullet(\tB)$ or $G^{\bullet}(\tB)$ to show that dependence on the extended exchange matrix. 

\subsection{Standard classes in terms of the Gysin spectral sequence} \label{sec:StandardGysin}
As summarized in Theorem~\ref{thm:LS}, the part of $H^k(\cA)$ in mixed Hodge degree $(k,k)$ is called the standard part of $H^{\ast}(\cA)$, and denoted $H^{\ast}(\cA)_{st}$.
%, and is generated as a ring by the classes $\dlog x_j$, where $x_j$ is a frozen variable, and by GSV-forms indexed by the connected components of $\Gamma$.
We now describe how to see the standard classes in $E_{\Gysin}$. 

As stated in Theorem~\ref{thm:LS}, $H^{\ast}(\cA)_{st}$ can be identified with the subring of differential forms in $\CC \langle \dlog x_1, \ldots, \dlog x_{n+m} \rangle$ which extend regularly to $\cA$. 
The ring  $\CC \langle \dlog x_1, \ldots, \dlog x_{n+m} \rangle$ is $G^{\emptyset} = \bigoplus_q E^{0q}_{1, \Gysin}$, and the condition that a form extends regularly to $\cA$ is the condition that it is in the kernel of the map $E^{0q}_{1, \Gysin} \to E^{1q}_{1, \Gysin}$. 
In other words, the standard forms are exactly those which end up in the first column of $E_{2, \Gysin}$. 

In particular, we consider the generators for the ring $H^{\ast}(\cA)_{st}$ from the second part of Theorem~\ref{thm:LS}.
For $x_j$ a frozen variable, we have $\dlog x_j = \theta(\{j\}, \emptyset)$.
For a connected component $\Delta$ of $\Gamma$, the corresponding GSV form is a linear combination of forms $\theta(\{i,j\}, \emptyset)$, where either $i$ and $j$ are vertices of $\Delta$, or else $i$ is a vertex of $\Delta$ and $j$ is frozen.

\subsection{Edge classes in terms of the Gysin spectral sequence} \label{sec:EdgeGysin}
Let $(a,b)$ be an edge of $\Gamma$. 
In Section~\ref{sec:acyclic}, we associated to $(a,b)$ a class $\epsilon_{ab}$ in $H^{3, (2,2)}(\cA)$. 
We originally defined $\epsilon_{ab}$ using a Mayer-Vietores sequence; let us now reconsider it in terms of a Gysin sequence.

We have $\cA = D_{\{ a \}} \sqcup D_{\{b \}} \sqcup \cA_{\{ a,b \}}$, where $D_{\{a \}}$ and $D_{\{ b \}}$ are the divisors where $x_a$ and $x_b$ vanish and $\cA_{\{ a,b \}}$ is their complement. 
We also recall  $\cA_{\{ a \}} = D_{\{ a \}} \sqcup \cA_{\{ a,b \}}$ and  $\cA_{\{ b \}} = D_{\{ b \}} \sqcup \cA_{\{ a,b \}}$. 
We have a Mayer-Vietores sequence coming from cover $\cA = \cA_{\{ a \}} \cup \cA_{\{ b \}}$ and a Gysin sequence coming from the decomposition $\cA = \cA_{\{ a,b \}} \cup \left( D_{\{ a \}} \sqcup D_{\{b\}} \right)$, and they fit together in a commutative diagram:
\[ \xymatrix@C=0.3 in{
\cdots H^k(\cA) \ar[r] \ar@{=}[d] & H^k(\cA_{\{ a \}}) \oplus H^k(\cA_{\{ b \}}) \ar[r] \ar[]+<-1.5em,-0.6em>;[d]^{\rho} & H^k(\cA_{\{ a,b \}}) \ar[r]^{\delta_{MV}} \ar[]+<0em,-0.6em>;[d]+<1.5em,0.6em>^{\sigma} & H^{k+1}(X) \ar@{=}[d]  \cdots\\
\cdots  H^k(\cA) \ar[r] & H^k(\cA_{\{ a,b \}})  \ar[r] & H^{k-1}(D_{\{a\}}) \oplus H^{k-1}(D_{\{b\}}) \ar[r]^-{\delta_{\Gysin}} & H^{k+1}(X)  \cdots\\
}\]
where the arrow labeled $\rho$ sends $(\alpha, \beta)$ to the restriction of $\alpha$ to $\cA_{\{a,b \}}$ and the arrow labeled $\sigma$ sends $\gamma$ to $(0, \mathrm{Residue}_{D_{\{b\}}}(\gamma))$. 
%(We hope that the angle of these arrows helps the reader recall which component of the direct sum they involve.)
Therefore, 
\[ \epsilon_{ab} = \delta_{MV}(\dlog x_a \wedge \dlog x_b) = \delta_{\Gysin}(0,  \mathrm{Residue}_{D_{\{b\}}}(\dlog x_a \wedge \dlog x_b)) = \delta_{\Gysin}(0, \dlog x_a) . \]

It is now straightforward to pass to the Gysin spectral sequence using more strata: $\epsilon_{ab}$ is represented by the cochain in $E^{12}_{1,\Gysin}$ corresponding to $\dlog x_a$ on the stratum $\cO_{\{b \}}$ or, in other words, the class $\theta(\{ a \}, \{ b \})$ in $G^{\{ b\}}$. 

\subsection{Rational acyclic exchange matrices} 
A \newword{(extended) rational exchange matrix} $M$ is a $(n+m) \times m$ matrix with rational entries such that the top $n \times n$ part of $M$ is skew-symmetric.  We observe that the complex $G^\bullet$ can be defined for a rational extended matrix $M$, replacing an extended exchange matrix $\tB$.  Write $G^\bullet(M)$ for the associated complex.
It will be useful to allow such matrices as intermediate steps, although they will not arise in our final results. 
This section collects lemmas about rational exchange matrices.

\begin{lemma}\label{lem:Mprime}
Let $M$ be an $(n+m) \times n$ rational exchange matrix and let  $M' =   UM$ where $U = \begin{sbm} \Id_n & 0 \\ P & Q \end{sbm}$
with  $P$ and $Q$ rational matrices of sizes $m \times n$ and $m \times m$, and $Q$ invertible.  Then the complexes $G^\bullet(M)$ and $G^\bullet(M')$ are naturally isomorphic.
\end{lemma}
\begin{proof}
%Right multiplication by $C$ simply scales the forms $\alpha_i$ by a scalar which clearly preserves the complex $\tG^\bullet(M)$. 
 Left multiplication by $U$ corresponds to a linear transformation \[ (\dlog x_1,\ldots, \dlog x_n, \dlog ,y_1,\ldots, \dlog y_m) \mapsto (\dlog x'_1,\ldots,\dlog x'_n, \dlog y'_1,\ldots,\dlog y'_m) \] where $\dlog \y$ and $\dlog \y'$ are related by an invertible linear transformation, and we have $\dlog x'_i -\dlog x_i \in {\rm span}(\dlog y_i)$.   It is straightforward to check that such a linear transformation commutes with the maps $\rho_I^J$ in the definition of the complex $G^\bullet(M)$.
\end{proof}

The following lemma is a variant of \cite[Proposition 5.11]{LS}.
Let $S^{\bullet}$ be the complex of graded vector spaces where the space $S^0$ is $\CC$ in degree $0$, the space $S^1$ is $\CC$ in degree $1$,  all other spaces are $0$ and all maps are $0$.

\begin{lemma} \label{lem:MRedundantRow}
Let $M$ be an $(n+m) \times n$ rational exchange matrix of full rank and let $M'$ be a  rational exchange matrix obtained by adding one more row to $M$.
Then $G^{\bullet}(M')$ is the tensor of $G^{\bullet}(M)$ with $S^{\bullet}$.
\end{lemma}

\begin{proof}
Since $M$ is of full rank, the additional row added to $M$ is in the span of the rows of $M$. So we can multiply $M'$ on the left by a matrix as in Lemma~\ref{lem:Mprime} to obtain the matrix $\begin{sbm} M \\ 0 \end{sbm}$. So $G^{\bullet}(M') \cong G^{\bullet}(\begin{sbm} M\\ 0 \end{sbm})$. Letting $x_{n+m+1}$ denote the cluster variable for the bottom row of $\begin{sbm} M \\ 0 \end{sbm}$, we see that $G^{\bullet}(\begin{sbm} M\\ 0 \end{sbm})$ is the tensor product of $G^{\bullet}(M)$ with the complex which has $1$ in degree $(0,0)$, and $\dlog x_{n+m+1}$ in degree $(1,1)$; this complex is isomorphic to $S^{\bullet}$.
\end{proof}

The cohomology of $G^{\bullet}(M)$ is essentially determined by the mutable part of $M$. 

\begin{lemma} \label{lem:MMutablePart}
Let $M_1$ and $M_2$ be $(n+m_1) \times n$ and $(n+m_2) \times n$ rational exchange matrices of full rank with the same mutable part. Let $m_1 \leq m_2$. Then $G^{\bullet}(M_2)$ is the tensor product of $G^{\bullet}(M_1)$ with $m_2-m_1$ copies of $S^{\bullet}$.
\end{lemma}

\begin{proof}
Let $M$ be the $(n+m_1+m_2) \times n$ matrix which has the same mutable part as $M_1$ and $M_2$, and has the frozen rows of $M_1$ stacked on top of the frozen rows of $M_2$.
Then we have $G^{\bullet}(M_1) \otimes (S^{\bullet})^{\otimes m_2} \cong G^{\bullet}(M) \cong G^{\bullet}(M_2) \otimes (S^{\bullet})^{\otimes m_1}$. 
However, tensoring with $S^{\bullet}$ is invertible on complexes supported in finitely many degrees, so we deduce that $G^{\bullet}(M_1) \otimes (S^{\bullet})^{\otimes (m_2-m_1)}$.
\end{proof}

Thus, given a skew-symmetric matrix $B$, we can append whatever rows to it we find most convenient and compute $G^{\bullet}$ for that choice of rows.
The following particular case will prove important in the next section. 
\begin{lemma} \label{lem:MIsolatedVertex}
Let $B$ be an $n \times n$ skew symmetric matrix whose $i$-th row and column are $0$. Let $B'$ be the $(n-1) \times (n-1)$ skew symmetric matrix obtained by deleting that zero and column from $B$. Let $M$ and $M'$ be rational exchange matrices of full rank obtained by adding $m$ and $m+1$ frozen rows, respectively, to $B$ and $B'$. Let $(G')^{\bullet}$ be the subcomplex of $G^{\bullet}(M)$ on those anticliques containing $i$. Then $(G')^{k+1, s+1} \cong G^{k,s}(M')$. 
\end{lemma}

\begin{proof}
Using Lemma~\ref{lem:MMutablePart}, we may assume that $M$ is of the form $\begin{sbm} M'& 0 \\ 0 &1 \end{sbm}$ where we have put the $0$-th column of $M$ in the last position. 
We write $x_{n+m}$ for the cluster variable corresponding to the last row of $M$.

The hypothesis on $B$ means that $i$ is an isolated vertex in the quiver $\Gamma$ of $M$, so, for every anticlique $I'$ for $M'$ extends to an anticlique $\{ i \} \cup I'$ for $M$, so the anticliques giving rise to the complex $G'$ are in bijection with the anticliques in the complex $G^{\bullet}(M')$. Thus, we must check that, for each anticlique $I'$ of $M'$, we have $G^{I \cup \{ i \},s+1}(M) \cong G^{I,s}(M')$, with the isomorphism commuting with the restriction maps $\rho^J_I$. Indeed, the isomorphism is given by $\wedge \dlog x_{n+m}$.
\end{proof}

\subsection{Reduction to the really full rank case}
\label{sec:decomposition}

As we have described, the complex $\tG^\bullet$ breaks into subcomplexes $\tG^\bullet[\chi]$ according to the elements $\chi \in X^{\ast}(I)$.
The goal of this section will be to reduce the computation of $\tG^\bullet(\tB)[\chi]$ to the computation of $G^{\bullet}(\tB')$ for another exchange matrix $\tB'$. 
We note that our proof will use rational exchange matrices as intermediate steps, but our main result (Theorem~\ref{thm:ReduceToG}) gives a matrix with integer entries.
Throughout this section, fix an exchange matrix $\tB$ and an element $\chi \in X^{\ast}$, and we write $B$ for the top part of $\tB$.

%The goal of this section will be to show that $\tG^\bullet(\tB)[g]$ is isomorphic to the complex $G^{\bullet}(\tB')$ for a generalized acyclic exchange matrix $\tB'$.
%Define $\tG^\bullet := G^\bullet[e]$ to be the summand of $G^\bullet$ associated to the trivial character $e \in X^*$.

We first note a case in which $\tG^\bullet[\chi]$ is quite trivial.
\begin{lem}
Let $\chi$ be an element of $X^{\ast}$ which is not contained in $X^{\ast}(I)$ for any anticlique $I$. Then $\tG^\bullet[\chi]$ is the zero complex, and the $\chi$-eigenspace of $H^{\ast}(\cA)$, for the action of $\Aut(A)$, is trivial.
\end{lem}

An example of a cluster variety for which the hypothesis of this lemma occurs is the cluster variety with exchange matrix $\begin{sbm} 0&2 \\ -2&0 \end{sbm}$: The group $X^{\ast}$ is $(\ZZ/2 \ZZ) \times (\ZZ/2 \ZZ)$; the anticliques are $\emptyset$, $\{ 1 \}$ and $\{ 2 \}$, and the corresponding subgroups are $\{ 0 \} \times \{ 0 \}$, $\{ 0 \} \times (\ZZ/2 \ZZ)$ and $(\ZZ/2 \ZZ) \times \{ 0 \}$, so the element $(1,1)$ is not in any of these groups.

\begin{proof}
In the formula $\tG^{ps} = \bigoplus_{I \in \cI_p} (G^I)^s \otimes \CC[X^{\ast}(I)]$, we only get a contribution to the degree $\chi$ if $\chi \in X^{\ast}(I)$. 
So, if $\chi$ is not contained in any $X^{\ast}(I)$, then $\tG^\bullet[\chi]=0$.
The $\chi$-eigenspace of $H^{\ast}(\cA)$ on $\Aut(A)$ is computed by the cohomology of $\tG^\bullet[\chi]$. 
\end{proof}

Now, suppose that $\chi$ is an element of $X^{\ast}$ which does lie in some $X^{\ast}(I)$. 
By Lemma~\ref{lem:MinimalJg}, there is a minimal anticlique $J(\chi)$ such that $\chi \in X^{\ast}(I)$ if and only if $I \supseteq J(\chi)$.
Define
\[ K(\chi) = \{ i \in [n] \setminus J(\chi) \mid \tB_{ij}=0 \ \text{for all}\ j \in J(\chi) \}. \]
%Our eventual goal will be to reduce the study of $\tG^\bullet[g]$ to the computation of $G^{\bullet}(M)$ for some matrix $M$ whose columns are indexed by $K(g)$. 
%To start with, we reduce to a matrix whose columns are indexed by $J(g) \cup K(g)$.
So $\tB$ has zeroes in the positions shown below, where we have ordered the rows and columns as $J(\chi)$, then $K(\chi)$,  then other mutable elements, then (in the case of rows) other mutable rows:
\[
\tB=
\begin{bmatrix}
0&0&\ast \\
0&B'&\ast \\
\ast&\ast&\ast \\
\hline
\ast & \ast & \ast \\
\end{bmatrix} 
\]
%\[
%\begin{bNiceMatrix}[\text{first-row}, \text{first-col}]
%& J(g) & K(g) & \text{other mutable} \\
%J(g) & 0 & 0 & \ast \\
%K(g) & 0 & \ast & \ast \\
%\text{other mutable} & \ast & \ast & \ast \\
% \hline
% \text{frozen} & \ast & \ast & \ast \\ 
% \end{bNiceMatrix}\]

Here is our main result:
\begin{theorem} \label{thm:ReduceToG}
With notation as above, let $k = |J(\chi)|$ and let $B'$ be the skew symmetric matrix obtained by restricting $B$ to the rows and columns indexed by $K(\chi)$. Let $\tB'$ be any full rank exchange matrix whose top part is $B'$ and has the same number of rows as $\tB$. Then the complex $\tG^{\bullet+k, s+k}[\chi]$ is isomorphic to $G^{\bullet,s}(\tB')$. 
\end{theorem}

\begin{proof}
As an intermediate step, let $\tB''$ be the  the matrix obtained from $\tB$ by taking the columns indexed by $J(\chi) \cup K(\chi)$, so the rows indexed by $[n] \setminus (J(\chi) \cup K(\chi))$ are relabeled as frozen. We get summands of $\tG^\bullet[\chi]$ from those anticliques $I$ with $I \supseteq J(\chi)$. All of these anticliques are contained in $J(\chi) \cup K(\chi)$. Choosing to delete a column indexed by some $i \not\in I$ does not change $G^I$. So we can delete the columns indexed by $[n] \setminus (J(\chi) \cup K(\chi))$ without changing the complex at all. 
Write $B''$ for the mutable part of $\tB''$.

We have thus reduced our attention to matrices of the form below, where the rows and columns are ordered with $J(\chi)$ first, then $K(\chi)$, and then (in the case of rows) the frozen rows:
\[
\tB''=
\begin{bmatrix}
0&0 \\
0&B' \\
\hline
\ast & \ast  \\
\end{bmatrix} 
\]
%\[ \tB''=
%\begin{bNiceMatrix}[\text{first-row}, \text{first-col}]
%& J(g) & K(g)  \\
%J(g) & 0 & 0  \\
%K(g) & 0 & \ast  \\
% \hline
% \text{frozen} & \ast & \ast \\ 
% \end{bNiceMatrix}\]
The degree $\chi$ part is the subcomplex of $G^{\bullet}(\tB'')$ on those anticliques containing $J(\chi)$.

We now apply Lemma~\ref{lem:MIsolatedVertex} to each of the columns indexed by $J(\chi)$. This shows that the subcomplex of $G^{\bullet}(\tB'')$ on those anticliques containing $J(\chi)$ is the complex $G^{\bullet}(\tB')$, with appropriate shifts of cohomological degree and grading, as required.
\end{proof}

Therefore, we will focus on computing $G^{\bullet}(\tB)$ in the rest of the paper, and we will use Lemma~\ref{lem:MMutablePart} to focus on the mutable part of $\tB$.

\section{Filtration of the Gysin complex}\label{sec:filtration}

We remind the reader of the notation $\theta(A, I)$ from Section~\ref{sec:basis}. So $G^I$ is spanned by $\{ \theta(A,I) \}_{A \subseteq [n+m] \setminus I}$. 
We now introduce a filtration of the complex $G^{\bullet}$ by the number of mutable elements in $A$: let $F^e G^I$ be the subspace of $G^I$ spanned by $\theta(A, I)$ for $|A \cap [n]| +|I| \geq e$.
We set $F^e G^{\bullet}$ to be the subspace of $G^{\bullet}$ spanned by the $F^e G^I$, as $I$ ranges over all anticliques.
Clearly, $F^0 G^{\bullet} \supseteq F^1 G^{\bullet} \supseteq \cdots$.

\begin{lemma} \label{lem:subcomplex}
The subspace $F^e G^{\bullet}$ is a subcomplex.
\end{lemma}

\begin{proof}
This is clear from the formula in Lemma~\ref{lem:rhoMapBasic}: we have $|(A \setminus \{ j \}) \cap [n]| = |A \cap [n]|-1$ and $|J| = |I|+1$, so $|(A \setminus \{ j \}) \cap [n]|+|J| = |A|+|I|$.
\end{proof}

As in Section~\ref{sec:basis}, let $N$ be a map from anticliques to subsets of $[n+m]$ such that $\tB_{N(I), I}$ is invertible for all anticliques $I$. 
We showed in Lemma~\ref{lem:Gbasis} that $\{ \theta(A, I) : A \subseteq I \sqcup N(I) \}$ is a basis of $G^I$. We now show that, if $N(I)$ is well chosen, this basis respects the filtration $F^{\bullet} G^I$ of $G^I$.

\begin{lemma} Let $I \in \cI$ be an anticlique.  Suppose that we have chosen $N(I)$ such that $|N(I)| \cap [n]$ is as small as possible, meaning that, for any subset $R$ of $[n+m]$ such that $\tB_{RI}$ is invertible, we have $|R \cap [n]| \geq |N(I) \cap [n]|$.
Then $\{ \theta(A, I) : |A \cap [n]|+|I| \geq e,\ A \subseteq [n+m] \setminus (I \sqcup N(I)) \}$ is a basis of $F^e G^I$. 
\end{lemma}

\begin{proof}
By Lemma~\ref{lem:Gbasis}, $\{ \theta(A, I) :  A \subseteq [n+m] \setminus (I \sqcup N(I)) \}$ is linearly independent, so the subset listed above is as well.
Clearly, all the listed basis elements are in $F^e G^I$. 
What remains is to show that these span $F^e G^I$. 
In other words, we need to show that, for any $A \subseteq [n+m] \setminus I$, we can write $\theta(A, I)$ as a linear combination of $\theta(A', I)$ with $A' \subseteq [n+m] \setminus (I \sqcup N(I))$ and $| A' \cap [n] | \geq |A \cap [n] |$.

Let $\tB_{iI}$ denote the $i$-th row of the matrix $\tB_I$. 
The vectors $\{ \tB_{iI} : i \in N(I) \}$ must be a basis of $\QQ^I$. 
Since we chose $[n] \cap N(I)$ as small as possible, we have that $[n+1, n+m] \cap N(I)$ is as large as possible, so $\tB_{N(I) \cap [n+1, n+m],\ I}$ must have the same rank as $\tB_{[n+1, n+m],\ I}$. 
In other words, for any $j \in [n+1, n+m]$, the row $\tB_{jI}$ must be in the span of $\{ \tB_{iI} : i \in N(I) \cap [n+1, n+m] \}$. 

Let $C = \tB_I (\tB_{N(I) I})^{-1}$, and let $C_j$ be the $j$-th row of $C$. The matrix $C$ has an identity matrix in positions $N(I) \times I$. 
Partition $I = I_1 \sqcup I_2$, where the $1$'s of this identity matrix in the columns indexed by $I_1$ lie in mutable rows and the $1$'s in the columns indexed by $I_2$ lie in frozen rows.
Since $C$ is obtained from $\tB_I$ by right multiplying by an invertible matrix, the relations between the rows of $C$ are the same as those between the rows of $\tB$. 
Thus, for any  $j \in [n+1, n+m]$, the row $C_j$ must be in the span of $\{ C_i : i \in N(I) \cap [n+1, n+m] \}$. 
Since $C$ has an identity matrix in positions $N(I) \times I$, this means that $C_{[n+1, n+m],I_1}$ must be $0$. 
Each column of $C$ gives a linear relation between the $\dlog x_i$, so we deduce that, if $i \in N(I) \cap [n]$, then $\dlog x_i$ is a linear combination of $\{ \dlog x_j : j \in [n] \setminus N(I) \}$. Meanwhile, if $i \in N(I) \cap [n+1, n+m]$, then  $\dlog x_i$ is a linear combination of $\{ \dlog x_j : j \in [n+m] \setminus N(I) \}$. 

We now return to the goal from the first paragraph: writing $\theta(A, I)$ as a linear combination of $\theta(A', I)$ with $A' \subseteq [n+m] \setminus (I \sqcup N(I))$ and $| A' \cap [n] | \geq |A \cap [n] |$. Indeed, we can use the relations from the previous paragraph to replace each $\dlog x_a$ for $a \in A \cap N(I) \cap [n]$ by a linear combination of $\{ \dlog x_j : j \in [n] \setminus N(I) \}$, and to replace $\dlog x_a$ for $a \in A \cap N(I) \cap [n+1, n+m]$ by a linear combination of $\{ \dlog x_j : j \in [n+m] \setminus N(I) \}$. Expanding the wedge of these replaced $1$-forms gives the required expression for $\theta(A,I)$.
\end{proof}

The filtered complex $F^\bullet G^\bullet$ simplifies significantly if $\tB$ has principal coefficients, as we now discuss.

\section{Principal coefficient case} \label{sec:principal}
\subsection{Reduction to principal coefficients} \label{PrincipalReduce}
Up to torus factors, any acyclic cluster variety of really full rank is isomorphic to the cluster variety with principal coefficients of the same exchange type.

\begin{proposition}[{\cite[Proposition 5.11]{LS}}]\label{prop:LSreduce}
Let $\cA(\tB)$ be an acyclic cluster variety of really full rank and let $B$ be the top part of extended exchange matrix $\tB$.  Then for some $a,b$ we have $\cA(\tB) \times (\GG_m)^a \cong \cA(B_\prin) \times (\GG_m)^b$.
\end{proposition} 

Using Proposition~\ref{prop:LSreduce}, we reduce the study of cluster algebras with really full rank to the study of cluster varieties with principal coefficients.  
We remind the reader that, in Theorem~\ref{thm:ReduceToG}, we have already reduced the study of full rank cluster algebras to the study of really full rank cluster algebras.

\subsection{Relation to cohomology of the independence complex}
We now assume that $\cA = \cA(B_\prin)$ is an acyclic cluster variety with principal coefficients.
Thus $m = n$, and we have cluster variables $x_1,\ldots,x_n$ and frozen variables $x_{n+1} = y_1,\ldots, x_{2n} = y_n$.   
%In Lemma \ref{lem:Sbasis}, we always have $r(I) = 0$, and we set $S(I) := I+n$.
For any subset $I$ of $[n]$, we take $N(I) = \{ i+n : i \in I \}$, so $\tB_{N(I),\ I}$ is an identity matrix. Thus, our basis of $G^I$ is indexed by subsets $A$ of $[2n] \setminus (I \sqcup N(I))$. 
It will be convenient to decompose such an $A$ as $C \sqcup (D+n)$ where $C = A \cap [n]$ and where $D = \{ d \in [n] : d+n \in A \}$. 
So $C$ and $D$ are subsets of $[n] \setminus I$.
Our basis of $G^I$ thus consists of the forms
\[ \theta(C,D,I):= \bigwedge_{c \in C} \dlog(x_c) \wedge \bigwedge_{d \in D} \dlog(y_d) \wedge \bigwedge_{i \in I} \alpha_i  \qquad \text{ where} \qquad C \cap I= D \cap I = \emptyset.\]
  Here, the factors are wedged together in the natural increasing order on $C$, $D$, and $I$.  

%Our basis for $G^I$ consists of the forms
%$$
%\theta(A,B,I):=\bigwedge_{a \in A} \dlog(x_a) \wedge \bigwedge_{b \in B} \dlog(y_b) \wedge \bigwedge_{i \in I} \alpha_i,$$ 

The differential $\partial$ of $G^\bullet$ sends
\begin{equation} \theta(C, D, I) \mapsto \sum_{\substack{c \in C \\ I \cup c \in \cI}}
\left.
 \begin{cases} \pm \theta(C {\setminus} c,\ D,\ I \cup c) & c \not \in D \\ 
 \sum_{\substack{d \in [n] {\setminus} C}} \pm \tB_{dc} \theta(C \cup d {\setminus} c,\ D {\setminus} c,\ I \cup c) & c \in D \\ \end{cases} \right\}. \label{Boundary} \end{equation}
The filtration $F^\bullet G^\bullet$ is by the cardinality of $C \cup I$. 

Let $\gr G^{\bullet} = \gr_F G^\bullet$ be the associated graded complex of the filtered complex $F^\bullet G^\bullet$.
The differential $\partial$ of $\gr G^{\bullet}$ is given by
\begin{equation} \partial(\theta(C, D, I)) = \sum_{\substack{c \in C {\setminus} D \\ I \cup c \in \cI}} \pm  \theta(C {\setminus} c, D, I \cup c) . \label{GradedBoundary} \end{equation}
We see from this equation that these maps preserve the values of $D$ and of $C \cup I$. Let $\gr G^{\bullet}(D, E)$ be the subcomplex where we fix $D$ and fix the value of $C \cup I$ to be $E$. Recall that $C \cap I = \emptyset$, so $E = C \sqcup I$.

Since the maps in $G^{\bullet}$ preserve the mixed Hodge structure, the maps in the associate graded complex also preserve the associated graded mixed Hodge structure.
Recall that $\theta(C,D,I)$ is in mixed Hodge degree $(s,s)$ for $s = |C|+|D|+|I| = |D|+|E|$. 
So the complex $\gr G^{\bullet}(D,E)$ is in mixed Hodge degree $(s,s)$ for $s = |D|+|E|$. 

The complex $\gr G^{\bullet}(D,E)$ has one term for each anticlique $I$ obeying $I \subseteq E {\setminus} D$. 
In other words, writing $\Gamma|_{E {\setminus} D}$ for the subgraph of $\Gamma$ induced on the vertices $E {\setminus} D$, the complex $\gr G^{\bullet}(D,E)$ computes the reduced cohomology of the simplicial complex $\cI(E {\setminus} D)$. 
We record the details below.

\begin{theorem}\label{thm:E1}
We have
\[ H^k(\gr G^{\bullet}(D,E)) \cong \tH^{k-1}(\cI(E {\setminus} D)) \]
where $\tH$ denotes reduced cohomology, and we use the convention that $\tH^{-1}(\{ \emptyset \})=\CC$ and $\tH^{-1}(\Delta)=0$ for every other simplicial complex.
\end{theorem}
%\begin{proof}
%For each independent set $J \in \cI_k(E {\setminus} D)$, we have a $(k-1)$-cochain $\psi_J$ taking the value one on the cycle $J$ and the value $0$ on all other cycles.  The (reduced) simplicial cochain complex $C^\bullet(\cI(E{\setminus} D))$ of $\cI(E{\setminus} D)$ has a basis given by the $\psi_J$, with differential
%$$
%\partial \psi_J = \sum_{\substack{J' \supset J \\ J' \in \cI_{k+1}}} (-1)^{J,J'} \psi_{J'}
%$$
%where the sign $(-1)^{J,J'}$ is equal to $(-1)^{k+1-r}$ if $J' = \{j_1 < j_2 < \cdots < j_{k+1}\}$ and $J = J' {\setminus} \{j_r\}$.  
%With $C$ fixed, define $(-1)^{A}$ to be $(-1)^{r_1+r_2+\cdots + r_{|A|}}$ where $A= \{c_{r_1},c_{r_2},\ldots,c_{r_|A|}\}$ and the $c$-s are ordered by $C = \{c_0> c_1> \cdots\}$.
%Checking the signs, we see that the map 
%\begin{equation}\label{eq:toIC}
%\Psi(\theta(A,B,I)) =(-1)^{A} (-1)^{|A||B|} \psi_I
%\end{equation}
%is an isomorphism of chain complexes $\Psi:  \gr G^{\bullet}(B,C) \to C^{\bullet-1}(\cI(C{\setminus} B))$, proving the proposition.
%\end{proof}
%We review some results on the topology of independence complexes in Section~\ref{sec:IC}.
%By Proposition~\ref{prop:E1}, if $\cI(\Gamma_{C {\setminus} B})$ is contractible, then $\gr G^{\bullet}(B,C)$ is exact. If $\Gamma_{C {\setminus} B}$ has an isolated vertex $v$, then $v$ will be a cone vertex of  $\cI(\Gamma_{C {\setminus} B})$ (see Corollary~\ref{cor:isolated}), so $\gr G^{\bullet}(B, C)$ will be exact whenever this occurs. 

We therefore pause to discuss results on the topology of independence complexes.

\subsection{Cohomology of independence complexes}\label{sec:IC}
Let $G$ be a finite undirected graph.  The independence complex $\cI(G)$ is the simplicial complex on the vertex set $V(G)$ whose faces are the independent sets in $G$. 

For a general graph $G$,  the independence complex $\cI(G)$ can have the homotopy type of any finite simplicial complex. 
Indeed, let $\Delta$ be a finite simplicial complex and let $\mathrm{Bary}(\Delta)$ be the first barycentric subdivision of $\Delta$. 
Let $G$ be the graph whose vertex set is the vertices of $\mathrm{Bary}(\Delta)$ and whose edges are the non-edges in the $1$-skeleton of  $\mathrm{Bary}(\Delta)$.
Then $\cI(G) \cong \mathrm{Bary}(\Delta)$, which is homotopy equivalent to $\Delta$.
Nonetheless, there are interesting results about the topology of $\cI(G)$, which we now describe.

Let $\Delta_1$ and $\Delta_2$ be two simplicial complexes on disjoint vertex sets.
The \newword{join} $\Delta_1 * \Delta_2$ has faces $F = F_1 \sqcup F_2$, where $F_i \in \Delta_i$.  

\begin{proposition}\label{prop:join}
Suppose $G$ is the disjoint union of $G_1$ and $G_2$.  Then 
\begin{equation}\label{eq:join}
\cI(G) = \cI(G_1) * \cI(G_2).
\end{equation}
\end{proposition}
Note that the join $S^a * S^b$ of two spheres is itself a sphere $S^{a+b+1}$.

\begin{corollary}\label{cor:isolated}
Suppose that $G$ contains an isolated vertex $v$.  Then $\cI(G)$ is contractible and $\tilde{H}^*(\cI(G))=0$.
\end{corollary}
\begin{proof}
In this case, $\cI(G)$ is the join of $\cI(G{\setminus} v)$ with a point, and is thus the cone over $\cI(G {\setminus} v)$.
\end{proof}

The following result is Theorem~7.1 in the preprint~\cite{MTarXiv}; it does not appear in the published version~\cite{MTpub}.

\begin{prop}\label{prop:vanishing} %(See Marietta and Testa, Theorem 7.1.). Disappeared from published version.
Let $G$ be a graph with $m$ vertices. Then $\tilde{H}^r(\cI(G))$ vanishes for $r>m/2-1$. 
\end{prop}

\begin{proof}
The proof is by induction on $m$.  The base cases $m=0$ and $m=1$ are clear. Suppose $m \geq 2$.  If $G$ has isolated vertices then by Corollary~\ref{cor:isolated}, we have $\tilde{H}^*(\cI(G))=0$, so the result holds.  If not, let $(u,v)$ be an edge of $G$. Let $G_u$, $G_v$ and $G_{uv}$ be the graphs where we restrict $G$ to the vertex sets $[m] {\setminus} u$, $[m] {\setminus} v$ and $[m] {\setminus} \{ u,v \}$. Then $\cI(G) = \cI(G_u) \cup \cI(G_v)$, and $\cI(G_{uv}) = \cI(G_u) \cap \cI(G_v)$. Writing out the Mayer-Vietores sequence, $\tilde{H}^r(\cI(G))$ sits between $\tilde{H}^{r-1}(G_{uv})$ and $\tilde{H}^r(G_u) \oplus \tilde{H}^r(G_v)$.  When $r > m/2-1$, all three of these cohomology groups vanish by induction, and therefore so does $\tilde{H}^r(\cI(G))$.
\end{proof}

We note that $\tilde{H}^{m/2}(\cI(G)) = \CC$ if $G$ is a union of $m/2$ disjoint edges,  in which case $\cI(G)$ is the boundary of an $m$-dimensional cross polytope, so this bound is optimal. Tracing through the proof, this is the only case where $\tilde{H}^{m/2}(\cI(G))$ is nonzero.
We now record some interesting results due to Ehrenborg and Hetyei and to Kozlov.

\begin{proposition}[{\cite[Corollary 6.1]{EH}}]\label{prop:EH1} 
Let $F$ be a forest.  Then $\cI(F)$ is either contractible or homotopy-equivalent to a sphere.
\end{proposition}

Let $F$ be a forest and $x \in V(F)$.  Then denote by $F_{x,h}$ the forest obtained by adding a path of length $h$ to the vertex $x$.  Thus $F_{x,0} = F$.  Also, let $F_{x,h,h'} := (F_{x,h})_{x,h'}$.  Then Ehrenborg and Hetyei show the following homotopy equivalences that recursively compute the homotopy type of independence complexes of forests.
\begin{proposition}[{\cite{EH}}] \label{prop:EH}
\begin{enumerate}
\item[(EH1)] $\cI(F_{x,1,1}) \simeq \cI(F_{x,1})$
\item[(EH2)] $\cI(F_{x,3}) \simeq \Sigma(F)$
\item[(EH3)] $\cI(F_{x,2,2}) \simeq \Sigma(F_{x,2})$
\item[(EH4)] $\cI(F_{x,2,1})$ is contractible.
\end{enumerate}
Here $\Sigma \Delta$ denotes the suspension of $\Delta$.
\end{proposition}
Thus, if $\Gamma$ is a forest, then every independence complex which arises in Theorem~\ref{thm:E1} will be contractible or a sphere. 
We mention what this result gives us for paths:

\begin{corollary}\label{cor:path}
Suppose that $G = P_n$ is a path of length $n$.  Then $\cI(P_n)$ is contractible if $n =3k$ and it is homotopy equivalent to $S^k$ if $n = 3k+1$ or $n = 3k+2$.
\end{corollary}

In contrast to forests, we have the following result for cycles.
\begin{proposition}[{\cite[Proposition 5.2]{Koz}}] \label{prop:Koz}
Let $C_m$ be an $m$-cycle.  Then the homotopy type of $\cI(C_m)$ is given by
$$
\cI(C_m) \simeq \begin{cases} S^{k-1} \vee S^{k-1} & \mbox{if $m = 3k$,}\\
S^{k-1} & \mbox{if $m= 3k\pm 1$.}
\end{cases}
$$ 
\end{proposition}

\subsection{Vanishing of mixed Hodge numbers}

We now prove Theorem~\ref{thm:main} stated in the introduction.
\begin{thm}\label{thm:vanishing}
Suppose $\cA$ is an acyclic cluster variety of really full rank.  Then $H^{k,(s,s)}(\cA) = 0$ for $s < (2/3)k$.
\end{thm}
%This is significantly stronger than the vanishing for $s < k/2$ coming from the fact that $\cA$ is smooth.
\begin{proof}
By Proposition~\ref{prop:LSreduce}, we reduce to the principal coefficient case.

Fix $s$.  According to Theorem~\ref{thm:explicitcomplex}, the $(s,s)$ part of $H^k(\cA)$ comes from $G^{\bullet,s}$.  The term $\theta(C, D, I)$ lives in degree $|C|+|D|+|I| = |D|+|E|$, so all contributions from $\gr G^{\bullet}(D,E)$ will live in mixed Hodge degree $(s,s)$ for $s=|D|+|E|$.  Now, $|E {\setminus} D| \leq |E| \leq |D|+ |E|$. So $(s,s)$ terms only come from $\gr G^{\bullet}(D,E)$ in cases where $|E {\setminus} D| \leq s$. 

By Proposition \ref{prop:vanishing}, $H^{r}(\gr G^{\bullet}(E,D))$ vanishes for $r>|E {\setminus} D|/2$.   It follows that the group $H^r(\gr G^{\bullet}(E,D))$ vanishes if $r>s/2$. Undoing the filtration, we see that $H^r(G^{\bullet, s})$ vanishes for $r>s/2$. Putting $r=k-s$, we deduce from Theorem~\ref{thm:explicitcomplex} that $H^{k, (s,s)}(\cA)=0$ for $k-s > s/2$ or, in other words, for $s<(2/3) k$. 
\end{proof}

\begin{remark}
In the boundary case $s=(2/3) k$, we have that $H^{k,(s,s)}(\cA)$ is spanned by the terms where $D = \emptyset$ and $\Gamma|_E$ is a disjoint union of $s/2$ edges.  We expect that $H^{3t, (2t,2t)}(\cA)$ is spanned by image of the the $t$-fold cup product map $H^{3, (2,2)}(\cA)^{\otimes t} \to H^{3t, (2t,2t)}(\cA)$.  We will show that $H^{3, (2,2)}(\cA) \cong H^1(\Gamma)$ in Proposition~\ref{prop:s=2}.
\end{remark}

\begin{remark}
A similar vanishing holds for $H^*(\cA)[\chi]$, for $\chi \in X^*$ and $\cA$ any full rank acyclic cluster variety.  In this case, the degrees are shifted by $|J(\chi)|$.  
%In the case $g = e$, we have $|J_g| = 0$, and we are in the situation of Theorem \ref{thm:vanishing}.
\end{remark}

\section{The spectral sequence of the filtered complex $F^\bullet G^\bullet$}\label{sec:E1}
\subsection{A spectral sequence}
We continue to assume that $\cA$ has principal coefficients.  The filtration $F^\bullet G^\bullet$ of the complex $G^\bullet$ gives a spectral sequence $E_{r, \Filt}^{ef}$ (not to be confused with the spectral sequence $E^{pq}_{r,\Gysin}$ in Section~\ref{sec:cohom}) that converges to the cohomology $H^{p+q}(G^\bullet)$. 
We remark that the spectral sequence $E^{ef}_{\Filt}$ is supported in the fourth quadrant, with $e \geq 0$ and $f \leq 0$. This matches standard conventions for the spectral sequence of a filtration; see for example \cite[\href{https://stacks.math.columbia.edu/tag/012K}{Tag 012K}]{stacks-project}.
We will often abbreviate $E_{r, \Filt}^{ef}$ to $E_r^{ef}$ for the rest of this paper, as $E_{\Gysin}$ will not reoccur.

All of the constructions we are about to make respect the splitting of $G^{\bullet}= \bigoplus_k G^{\bullet,s}$, and its associated objects, into mixed Hodge degrees.
We will write  $E_r^{ef} = \bigoplus_s E_r^{ef,s}$.

We recall the standard notations for the spectral sequence associated to a filtration $F^{\bullet} G^{\bullet}$ on a complex $G^{\bullet}$
In the $r$-th page $E_{r}$, the differentials have degree $(r,1-r)$.  The $0$-th page is 
\[
E_{0}^{ef} = \gr^e G^{e+f}: = F^e G^{e+f} / F^{e+1} G^{e+f}
\]
with vertical differentials those of the associated graded complex
\[
\gr^e G^\bullet := \gr^e G^0 \to \gr^e G^1 \to \cdots.
\]
So the $E_1$-page is given by
\[
E_{1}^{ef} = H^{e+f}(\gr^e G^\bullet)
\]
with horizontal differentials $\partial_1:H^{e+f}(\gr^e G^\bullet) \to H^{e+f+1}(\gr^{e+1}G^\bullet)$.  In general, the $E_r$-page is given by
\[
E_r^{ef} = \frac{\{x \in F^e G^{e+f} \mid \partial x \in F^{e+r} G^{e+f+1}\}}{F^{e+1}G^{e+f} + \partial(F^{e-r+1} G^{e+f-1})}.
\]
 The differential $\partial_r: E_r^{ef} \to E_r^{e+r,f+1-r}$ is simply that induced by $\partial:G^\bullet \to G^\bullet$.

\subsection{The $E_0$-page and the spaces on the $E_1$-page}
Spelling out our definitions, $E_0$ sits in the cone $\{ (e,f)  : 0 \leq e+f \leq e \}$.
A basis for $E_0^{ef}$ is indexed by triples $(C,D,I)$ where $|I| = e+f$, $|C| = -f$ and $C \cap I = D \cap I = \emptyset$ (so $|E| = |C|+|I| = e$). 
The differential maps are in the vertical direction $(0,1)$ and, if we restrict ourselves to a vertical line where $e$ is fixed, we see the cochain complexes for the reduced cohomology of the simplicial complexes $\cI(E \setminus D)$, where $|E| = e$. 
The degree $-1$ part of the cochain complex (which is always $\CC$, since we are using reduced cohomology) is on the line $e+f=0$.

The mixed Hodge degree, $s$, is $|D|+|E| = e+|D|$. 
The cohomology group $H^k(\cA)$ will involve terms with $k = |C|+|D|+2|I| =2e+f+|D|$.
We depict $E_{0,\Filt}$ in Figure~\ref{E0FiltFigure}.

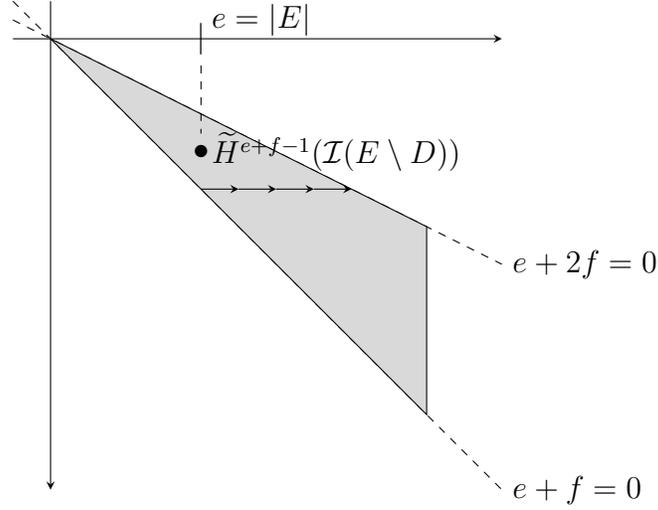
\begin{figure}
\begin{tikzpicture}
\draw[-stealth] (-0.5,0) -- (6,0) ;
\draw[-stealth] (0,0.5) -- (0,-6) ;
\draw[dashed] (-0.5,0.5) -- (6,-6) ;
\node[anchor=west] (diag) at (6,-6) {$e+f=0$};
\draw[dashed] (-0.5,0.25) -- (6,-3) ;
\node[anchor=west] (halfdiag) at (6,-3) {$e+2f=0$};
%\draw[dashed] (0,5) -- (5,0) ;
\draw[fill=gray!30] (0,0) -- (5,-2.5) -- (5,-5) -- (0,0);
\draw[-{stealth[scale=2]}] (2,-2) -- (2.5,-2) ;
\draw[-{stealth[scale=2]}] (2.5,-2) -- (3,-2) ;
\draw[-{stealth[scale=2]}] (3,-2) -- (3.5,-2) ;
\draw[-{stealth[scale=2]}] (3.5,-2) -- (4,-2) ;
\node (point) at (2,-1.5) {$\bullet$};
\node[anchor=west] (Htilde) at (2,-1.5) {$\widetilde{H}^{e+f-1}(\cI(E \setminus D))$};
\node (eTick) at (2,0) {$|$} ;
\node[anchor=west] (eLabel) at (2,0.25) {$e=|E|$} ;
\draw[dashed] (eTick) -- (point);
\end{tikzpicture}
\caption{The page $E_{1, \Filt}$. We have $s=e+|D|$ and $k=2e+f+|D|$. \label{E1FiltFigure}}
\end{figure}

The spaces on the $E_1$ page are the cohomology of the complexes on the $E_0$ page, so 
\[
E_1^{ef} = \bigoplus_{\substack{(D,E) \\ |E|=e}} \tH^{e+f-1}(\cI(E {\setminus} D)),
\]
and the Deligne splitting is 
\[
E_1^{ef,s} = \bigoplus_{\substack{(D,E) \\ |E|=e, \; |D|=s-e}} \tH^{e+f-1}(\cI(E {\setminus} D)) .
\]
From Proposition~\ref{prop:vanishing}, the page $E_{1, \Filt}$ is supported on $e/2 \leq -f \leq e$. 
We depict  $E_{1,\Filt}$ in Figure~\ref{E1FiltFigure}.

\begin{figure}
\begin{tikzpicture}
\draw[-stealth] (-0.5,0) -- (6,0) ;
\draw[-stealth] (0,0.5) -- (0,-6) ;
\draw[dashed] (-0.5,0.5) -- (6,-6) ;
%\draw[dashed] (0,5) -- (5,0) ;
\draw[fill=gray!30] (0,0) -- (5,0) -- (5,-5) -- (0,0);
\node (point) at (3,-1.5) {$\bullet$};
\node[anchor=west] (theta) at (3,-1.5) {$\theta(C,D,I)$};
\node (eTick) at (3,0) {$|$} ;
\node[anchor=west] (eLabel) at (3,0.25) {$e=|C|+|I|=|E|$} ;
\draw[dashed] (eTick) -- (point);
\node (fTick) at (0,-1.5) {$-$};
\node[anchor=east] (fLabel) at (0,-1.5) {$f=-|C|$} ;
\draw[dashed] (fTick) -- (point);
\node (ITick) at (1.5,0) {$|$} ;
\draw[dashed] (ITick) -- (point);
\node[anchor=south] (ILabel) at (1.5,0) {$e+f=|I|$} ;
\node (0point) at (3,-3) {$\bullet$};
\node[anchor=west] (0theta) at (3,-3) {$\theta(E,D,\emptyset)$};
\draw[-{stealth[scale=2]}] (3,-3) -- (3,-2.5) ;
\draw[-{stealth[scale=2]}] (3,-2) -- (3,-1.6) ;
\draw[-{stealth[scale=2]}] (3,-1.35) -- (3,-1) ;
%\node (qTick) at (0,2.5) {$\bullet$} ;
%\node[anchor=east] (qLabel) at (0,2.5) {$q=s=|A|+|I|$} ;
%\node (point) at (1,2.5) {$\bullet$};
%\node[anchor=west] (theta) at (1,2.5) {$\theta(A,I)$};
%\draw[dashed] (point) -- (pTick);
%\draw[dashed] (point) -- (qTick);
%\node[anchor=west] (dimEqn) at (5,0.25) {$p+q=m+n$};
%\node[anchor=west] (diagonal) at (5,5.25) {$p=q$};
%\draw[-{stealth[scale=2]}] (4.25,1.5) -- (4.75,1.5) ;
%\node[anchor=west] (map) at (4.75, 1.5) {\tiny the $\partial$-map};
\end{tikzpicture}
\caption{The page $E_{0, \Filt}$. We have $s=e+|D|$ and $k=2e+f+|D|$. \label{E0FiltFigure}}
\end{figure}

\subsection{Differentials on the $E_1$-page}
%By Proposition \ref{prop:E1}, we have an explicit description of all the groups $E_1^{pq}$ for the spectral sequence associated to the filtered complex $F^\bullet G^\bullet$.  Namely,
%\[
%E_1^{pq} = H^{p+q}(\gr^p G^\bullet) = \bigoplus_{\substack{(B,C) \\ |C|=p}} \tH^{p+q-1}(\cI(C {\setminus} B))
%\]
%%where the sum is over $(B,C)$ satisfying $|C| = p$.  
%This decomposes into the direct sum
%\begin{equation}\label{eq:E1}
%E_1^{pq,s} = \bigoplus_{\substack{ (B,C) \\ |B|=s-p \\ |C|=p }} \tH^{p+q-1}(\cI(C {\setminus} B))
%\end{equation}
%%where the sum is over $(B,C)$ satisfying $|B| + |C| = s$ and $|C| = p$.
%We note that the simplicial complex $\cI(C \setminus B)$ has $|C \setminus B| \leq |C| =p$ vertices, so $\tH^{\ast}(\cI(C \setminus B))$ is supported in degrees at most $p-1$, forcing $q \leq 0$.

We now want to describe the maps $\partial_1: E_1^{e,f,s} \to E_1^{e+1,f,s}$ on the $E_1$ page.  This differential $\partial_1$ is the sum of maps
\begin{equation}\label{eq:BC}
\partial_{(D,E) \to (D',E')}:
\tH^{e+f-1}(\cI(E {\setminus} D)) \to \tH^{e+f}(\cI(E' {\setminus} D')),
\end{equation} where $|D'| = |D|-1$ and $|E'| = |E|+1$.  
%Also recall that the cohomological degree here corresponds to $|I|-1$.
Here is our main result:
\begin{theorem} \label{thm:E1map}
If the map $\partial_{(D,E) \to (D',E')}$ is nonzero, then there are elements $a$ and $b$ of $[n]$, with $a \in D \cap E$, $b \in [n] \setminus (D \cup E)$ and $B_{ab} \neq 0$, such that $D' = D \setminus \{ a \}$ and $E' = E \cup \{b \}$. In that case, we will have $\cI(E' \setminus D') = \cI(E \setminus D') \cup \cI(E' \setminus D)$ and $\cI(E \setminus D) = \cI(E \setminus D') \cap \cI(E' \setminus D)$. The map $\partial_{(D,E) \to (D',E')}: \tH^{k}(\cI(E {\setminus} D)) \to \tH^{k+1}(\cI(E' {\setminus} D'))$ is $\pm B_{ab}$ times the boundary map in the Mayer-Vietores sequence coming from the cover $\cI(E' \setminus D') = \cI(E \setminus D') \cup \cI(E' \setminus D)$.
\end{theorem}

%\begin{remark}
%It follows from Theorem~\ref{thm:E1map} that $E_{1, \Filt}$ can be decomposed as a direct sum of subcomplexes where the set $E \setminus D$ is constant on each summand, as are the quantities $\#(C \setminus B) + 2 \#(B \cap C)$ and $\#(C \setminus B) + 2 \#([n] \setminus (B \cap C))$.
%\end{remark}

\begin{figure}
\begin{tikzpicture}
\draw (0,0) rectangle (7,5) ;
\draw[rounded corners=30 pt] (1,1) rectangle (4,4) ;
\draw[rounded corners=30 pt] (3,1) rectangle (6,4) ;
\draw[dashed] (4,3) arc (90:270:0.5);
\draw[dashed] (6,2) arc (270:450:0.5);
\node (D) at (2.5,3.5) {$D$};
\node (E) at (4.5,3.5) {$E$};
\node (a) at (3.75,2.5)  {$a$};
\node (b) at (6.25,2.5)  {$b$};
\end{tikzpicture}
\caption{A Venn diagam depicting $D$, $E$, $a$ and $b$. \label{Venn}}
\end{figure}
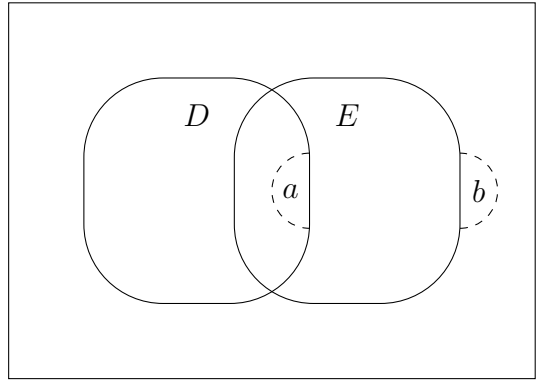

We depict $D$, $E$, $a$ and $b$ in a Venn diagram in Figure~\ref{Venn}. 
The boundaries of $D$ and $E$ are the rounded rectangles. 
To form $D'$ and $E'$, modify the boundaries of $D$ and $E$ to follow the dashed semicircles.

The proof of Theorem~\ref{thm:E1map} will occupy the rest of this section.

Let $\theta(C,D,I) \in F^e G^\bullet$.  Define 
\begin{equation}\label{eq:partiala}
\partial^{(a)} \theta(C,D,I) = \begin{cases}   \sum\limits_{b \in [n] \setminus (C \cup I)} \pm \tB_{ba}  \theta(C \cup b {\setminus} a, D {\setminus} a, I \cup a) &\mbox{if $a \in C \cap D$ and $I \cup a \in \cI$} \\
\qquad 0 & \mbox{otherwise.}
\end{cases}
\end{equation}
The maps $\partial^{(a)}$ induce the differentials \eqref{eq:BC}.  Namely, for $\theta \in E_1^{e,\bullet} = F^e G^\bullet$, we have
\[
\partial \theta =  \sum_{a} \partial^{(a)} \theta  \in F^{e+1} G^\bullet.
\]

% By \eqref{Boundary}, the projection of $\partial \theta(C,D,I)$ to $F^{e+1} G^\bullet$ is equal to
%$
%\sum_{a \in C \cap D} \partial^{(a)} \theta(C,D,I) 
%$
%where for $a \in C \cap D$, we have
%\begin{equation}\label{eq:partiala}
%\partial^{(a)} \theta(C,D,I) = \begin{cases} -  \sum\limits_{b \in [n] \setminus (C \cup I)} \tB_{ba} (\pm \theta(C \cup b {\setminus} a, D {\setminus} a, I \cup a)) &\mbox{if $I \cup a \in \cI$} \\
%\quad 0 & \mbox{otherwise.}
%\end{cases}
%\end{equation}
%These maps induce the differentials \eqref{eq:BC}.
%At this point, we can prove part of Theorem~\ref{thm:E1map}.

\begin{lemma} \label{lem:partialVanish}
If the map $\partial_{(D,E) \to (D',E')}$ is nonzero, then there is some $a \in D \cap E$ and $b \in [n] \setminus E$, such that such that $D' = D \setminus \{ a \}$, $E' = E \cup \{ b \}$ and $B_{ab} \neq 0$. 
\end{lemma}

\begin{proof}
Suppose that $\partial_{(D,E) \to (D',E')}$ is nonzero. Then there is some $\theta(C,D,I)$ with $E = C \sqcup I$ such that there is a nonzero summand $B_{ba} \theta(C \cup b {\setminus} a, D {\setminus} a, I \cup a)$ in $\partial \theta(C,D,I)$ with $D' = D \setminus a$ and $E' = C \cup I \sqcup \{b\} = E \sqcup \{ b \}$. 
Thus, we have deduced that $D' = D \setminus \{ a\}$, that $a \in D \cap E$ and that $E' = E \sqcup \{ b \}$. 
In addition, we have $B_{ba} \neq 0$ as required.
\end{proof}

What remains is (1) to show that $\partial_{(D,E) \to (D',E')}$ vanishes if $b \in D$ and, (2) in the case that $b \not\in D$, to show that  $\partial_{(D,E) \to (D',E')}$ is related to the Mayer-Vietores sequence as claimed. 
We tackle the first task now:

\begin{prop}\label{prop:addvertex}
Suppose that $a \in D \cap E$ and $b \in D \setminus E$. Put $D' = D \setminus \{ a \}$ and $E' = E \cup \{ b \}$. 
Then $\partial_{(D,E) \to (D',E')}=0$.
\end{prop}

\begin{proof}
It is convenient to set  $X = E {\setminus} D$ and $X' = E' {\setminus} D'$, so $X' = X \sqcup \{ a \}$.
Thus we are discussing a map from $\widetilde{H}^{k}(\cI(X))$ to $\widetilde{H}^{k+1}(\cI(X'))$.

Let $\psi \in C^k(\cI(X))$ represent a class in $H^k(\cI(X))$.
%, the image of some element of $G^\bullet$ under the map \eqref{eq:toIC}.
We may view $\psi$ also as a cochain on $\cI(X')$, which we denote
$\psi' \in C^k(\cI(X'))$.    Let $\psi_J$ denote the cochain taking the value one $J$ and the value $0$ on all other $J' \neq J$.
By assumption, $\partial \psi = 0 \in
C^{k+1}(\cI(X))$, so if $\psi = \sum_I z_I \psi_I$ then we must have
\[
\partial \psi' = \sum_{I \cap N[a] = \emptyset} \pm z_I \psi_{I \cup \{a\}}.
\]
The cochain $\partial \psi'$ is tautologically $0$ in
$H^{k+1}(I(X'))$.  Up to an overall factor of $\pm \tilde B_{ba}$,
this is the image of $\psi \in H^k(I(X))$ under the map
$\partial_{(D,E) \to (D',E')}$.  Thus $\partial_{(D,E) \to (D',E')}$ is 
$0$.
\end{proof}

Thus, the only nonzero maps are in the case $D' = D \setminus \{ a \}$ and $E' = E \cup \{ b \}$ for $a \in D \cap E$ and $b \in [n] \setminus (D \cup E)$, and $B_{ba} \neq 0$.
We have thus established all the vanishing claims from Theorem~\ref{thm:E1map}.
It remains to relate $\partial_{(D,E) \to (D',E')}$ to the Mayer-Vietores map.

We abbreviate $X := E \setminus D$. So $E \setminus D' = X \sqcup \{ a \}$, $E' \setminus D = X \sqcup \{ b \}$ and $E' \setminus D' = X \sqcup \{ a,b \}$.
%So we have containments of independence complexes \[ \cI(X) \ \subset \ \cI(X \sqcup \{ a \}), \cI(X \sqcup \{ b \}) \ \subset \ \cI(X \sqcup \{ a,b \}) \]
%with $\cI(X) =  \cI(X \sqcup \{ a \}) \ \cap \ \cI(X \sqcup \{ b \})$.
We note that the condition $B_{ba} \neq 0$ is equivalent to saying that $(a,b)$ is an edge of $\Gamma$, so no independence set in $E' \setminus D'$ contains both $a$ and $b$. 
So we have
\[ \cI(X \sqcup \{ a,b \}) =   \cI(X \sqcup \{ a \}) \ \cup \ \cI(X \sqcup \{ b \}) \quad \text{and} \quad \cI(X) =   \cI(X \sqcup \{ a \}) \ \cap \ \cI(X \sqcup \{ b \})  . \]
Thus, the Mayer-Vietores sequence in Theorem~\ref{thm:E1map} makes sense, and it remains to compare the map $\partial_{(D,E) \to (D',E')}$ to the boundary map $\delta$ from the Mayer-Vietores sequence.

The Mayer-Vietores map $\tH^k(\cI(X)) \to \tH^{k+1}(\cI(X \sqcup \{ a,b \}))$ is defined as follows.
Take an element of $\tH^k(\cI(X))$ and represent it by a cocycle $\psi \in C^k(\cI(X))$; in other words, $\psi$ is a function on $(k+1)$-element anticliques of $D$.
We define $\eta \in C^{k+1}(\cI(X \sqcup \{ a,b \}))$ roughly as follows:
\[ \eta(I) = \begin{cases} 
\psi(I \setminus \{ a \}) & a \in I \\
0 & a \not\in I \end{cases} .\]
We say ```roughly" because a cocycle is a function on ordered anticliques, which is anti-symmetric with respect to reordering. So, more precisely, we mean that these formulas are correct if $a$ is the last element of $I$ and $I \setminus \{ a \}$  is ordered by the induced order from $I$.
Then $\eta$ is a cocycle in $C^{k+1}(\cI(X \cup \{ a,b \}))$ and represents the image of the Mayer-Vietores map. 

We now check that this differs from $\partial_{(D,E) \to (D', E')}$ by a factor of $B_{ba}$. 
For a cocycle $\psi$, we write $[\psi, D,E]$ for the form $\sum_{I \in \cI(D)} \psi(I) \theta(E \setminus I, D, I)$ in $G^{\bullet}$. 
So
\[ \partial_{(D,E) \to (D', E')} [\psi, D, E] = \sum_{I \in \cI(X) : I \cup a \in \cI(X \cup a)} B_{ba} \psi(I) \theta{\big (}E \cup b \setminus (I \cup a), D \setminus a, I \cup a {\big )} . \]
We see that the right hand side is $B_{ba} [ \eta, D \setminus a, E \cup  b ]$ for the $\eta$ described above. This concludes the proof of Theorem~\ref{thm:E1map}.

\subsection{Frozen classes in the filtration spectral sequence}
In Section~\ref{sec:StandardGysin}, we saw that the standard classes in $H^{\ast}(\cA)$ came from the $I=\emptyset$ portion of $E_{1, \Gysin}$. 
%We now restate this in terms of the filtration spectral sequence.
The class $\theta(C, D, \emptyset)$ always lies on the antidiagonal $e+f=0$ of $E_{\Filt}$, so all the standard classes will lie on this antidiagonal.

\begin{prop} \label{E00}
For all $r \geq 0$, the space $E^{00}_{r, \Filt}$ is $2^n$ dimensional, with basis the forms $\theta(\emptyset, D, \emptyset)$. The corresponding subspace of $H^{\ast}(\cA)$ is the free exterior algebra on the $d \log y_i$.
\end{prop}

\begin{proof}
The space $E^{00}_{0, \Filt}$ has basis $\theta(\emptyset, D, \emptyset)$, for $D$ ranging over subsets of $[n]$. 
Each of these corresponds to the form $\bigwedge_{j \in D} \dlog y_j$ on the big cluster torus. 
Every such form maps to $0$ in $G^{\{ i \}}$ for each singleton $\{ i \}$, so in particular it maps to $0$ in $F^r G^{\{ i \}}$, and thus all maps out of the $(0,0)$ position are $0$ on the entire spectral sequence.
\end{proof}

%For each frozen variable $x_{i+n}$, the form $\dlog x_{i+n}$ is $\theta(\emptyset, \{ i \}, \emptyset)$. 
%This sits in $E^{00}_{0, \Filt}$ and survives to $E^{00}_{r,\Filt}$ for all larger $r$.
%Indeed, $E^{00}_{r,\Filt}$ is precisely the exterior algebra on the classes $\dlog y_i$.
\subsection{GSV classes in the filtration spectral sequence}
Let $\Delta$ be a connected component of $\Gamma$. The GSV form $\gamma_{\Delta}$ is $\sum_{i \in \Delta} \dlog x_i \wedge \dlog y_{i} + \sum_{i,j \in \Delta} \tB_{i,j} \dlog x_i \wedge \dlog x_j$, where we abuse notation by also writing $\Delta$ for the set of vertices of $\Delta$. The terms in the first sum are in $F^1 G^{\emptyset}$ and the terms in the second sum are in the smaller subspace $F^2 G^{\emptyset}$. So this corresponds to the term $\sum_{i \in \Delta} \theta(\{i \}, \{i \}, \emptyset)$ in $E^{1(-1)}_{0, \Filt}$.
This form survives to all later pages $E^{1(-1)}_{r, \Filt}$, and so the GSV forms span a subspace of $E^{1(-1)}_{r, \Filt}$ naturally isomorphic to $H^0(\Gamma)$. 
%Indeed, as we will see in REFERENCE, for $r \geq 2$, we have $E^{1(-1)}_{r, \Filt} = H^0(\Gamma) \wedge \Alt^{\bullet} \langle \dlog x_{n+i} \rangle$.  I think this sentence needs adjustment: if $|\Delta|=1$ then $\gamma = \dlog x_1 \wedge \dlog y_1$ so $\gamma \wedge \dlog y_1 = 0$.

\subsection{Edge classes in the filtration spectral sequence}
Let $(a,b)$ be an edge of $\Gamma$. In Section~\ref{sec:EdgeGysin}, we saw that the edge class $\epsilon_{ab} \in H^{3, (2,2)}(\cA)$ was represented by $\theta(A,I)$ with $A = \{ a \}$ and $I = \{ b \}$. 
In our current notation for the principal coefficients case, this is called $\theta(\{ a \}, \emptyset, \{ b \})$. Thus, it arises from $\widetilde{H}^0(\cI(\{ a,b \}))$ in $E^{2(-1)}_{1,\Filt}$. Since $(a,b)$ is an edge of $\Gamma$, the independence complex $\cI(\{ a,b \})$ is two points, and $\widetilde{H}^0(\cI(\{ a,b \}))$ is one dimensional, spanned by the edge class $\epsilon_{ab}$.

\section{Small mixed Hodge degrees}\label{sec:small}
In this section, let $\cA$ be an acyclic cluster variety of really full rank.  We calculate $H^{\ast,(s,s)}(\cA)$ for small values of $s$.  We use Theorem~\ref{thm:explicitcomplex} and carry out our calculations in the case of principal coefficients, and use Proposition~\ref{prop:LSreduce} to extend to the acyclic really full rank case.

Choose, once and for all, a GSV-form $\gamma_{\Gamma_i}$ for each connected component $\Gamma_i$ of $\Gamma$, and let $H_{\GSV} \subseteq H^{2, (2,2)}(\cA)$ be the vector space they span.  By the analysis in \cite[Section 9]{LS}, the $\gamma_{\Gamma_i}$ are linearly independent, and we have an isomorphism $H_{\GSV} \cong H^0(\Gamma)$.  Also, let $H_{\frozen} :=  \bigoplus_{i=1}^m \CC \dlog y_i \subset H^{1,(1,1)}(\cA)$.  

\begin{remark}
Recall from Theorem~\ref{thm:LS}(2) that the top weight cohomology $H^{\ast}(\cA)_{st}$ is generated by $H_{GSV}$ and $H_{\frozen}$.   It would be interesting to study $\bigoplus_k H^{k, (k-r, k-r)}(\cA)$ as a $H^{\ast}(\cA)_{st}$-module; 
unfortunately, the methods of this paper do not immediately yield this.
\end{remark}
%
%\begin{remark} According to Theorem~\ref{thm:vanishing}, $H^{\ast,(s,s)}(\cA) = H^{s,(s,s)}(\cA)$ when $s = 0,1$.  In \cite{LS}, the entire top weight cohomology $H^{s,(s,s)}(\cA)$ was computed.  
%\end{remark}
%
%Set $H_{\frozen} = H^{1, (1,1)} = H^1$. As we will verify in Proposition~\ref{prop:s=1}, a basis for $H_{\frozen}$ is $\dlog y_i$, where $y_1$, $y_2$, \dots, $y_m$ are the frozen variables. As we will now explain, the frozen terms account for the $(0,0)$ position in all pages of $E_{\Gysin}$:
%We begin by computing $E^{00}_{r, \Filt}$, as it will contribute to every other computation we perform.

%We now continue our computation, computing $H^{\ast, (s,s)}$ for the first few values of $s$:
%
%Besides the frozen forms, another easy to understand part of $H^{\ast}(\cA)$ are the Gehktman-Shapiro-Vainshtein forms.
%Recall that, for a connected component $\Gamma_i$ of $\Gamma$, we defined a GSV-form $\gamma_{\Gamma_i}$ in Section~\ref{sec:cohom}. 
%Choose, once and for all, a GSV-form $\gamma_{\Gamma_i}$ for each connected component of $\Gamma$.
%By the analysis in Section~9 of~\cite{LS}, the $\gamma_{\Gamma_i}$ are linearly independent; let $H_{\GSV} \subseteq H^{2, (2,2)}$ be the vector space they span. 
%
%In~\cite{LS}, the authors showed that $\bigoplus_k H^{k, (k,k)}(\cA)$ is generated as a ring by $H_{\frozen}$ and $H_{\GSV}$; we will likewise see in this section that other parts of the mixed Hodge structure of $\cA$ are best expressed in terms of $H_{\frozen}$ and $H_{\GSV}$.

\subsection{The case $s=0$} 
\begin{proposition}\label{prop:s=0}
Let $\cA$ be an acyclic cluster variety of really full rank.
Then $H^{\ast,(0,0)}(\cA) = H^{0,(0,0)}(\cA) =\CC$.
\end{proposition} 
The complex $G^{\bullet, 0}$ is given by 
$$
\CC \cdot \theta(\emptyset, \emptyset, \emptyset) \to 0 \to 0 \to \cdots
$$
Thus $H^{\ast,(0,0)}(\cA) = H^{0,(0,0)}(\cA) \cong \CC$.

\subsection{The case $s=1$}  
\begin{proposition}\label{prop:s=1}
Let $\cA$ be an acyclic cluster variety of really full rank.
Then we have $H^{\ast,(1,1)}(\cA) = H^{1,(1,1)}(\cA) =H_{\frozen} =  \bigoplus_{i=1}^m \CC \dlog y_i$.
\end{proposition} 

\begin{proof}
We assume that we are in the principal coefficient case.  Our complex is built out of the terms $\theta(C,D,I)$ with $|C| + |D| + |I|= |D| +|E| = 1$.   Now if $(|D|, |E|) = (0,1)$, then $\Gamma_{E {\setminus} D}$ is an isolated vertex, so the complex $G^{\bullet,1}_{gr}(D,E)$ is exact and does not contribute past page $E_0$. 
So the only contribution is from terms where $C = I = \emptyset$. 
These terms are in $E^{00}$, and we have already computed in Proposition~\ref{E00} that the part of $E^{00}$ in mixed Hodge degree $(1,1)$ is $\bigoplus_{i=1}^n \CC \dlog y_i$.
Proposition~\ref{prop:LSreduce} then lets us transfer this result to any other acyclic really full rank case.
\end{proof}

\subsection{The case $s=2$} 
\begin{proposition}\label{prop:s=2}
Let $\cA$ be an acyclic cluster variety of really full rank.  Then $H^{\ast,(1,1)}(\cA) = H^{2,(2,2)}(\cA) \oplus H^{3,(2,2)}(\cA)$.  We have $H^{2,(2,2)}(\cA) \cong \bigwedge^2 H_{\frozen} \oplus H_{\GSV}$, and $H^{3,(2,2)}(\cA) \cong H^1(\Gamma)$ is spanned by the edge classes.
\end{proposition}
Using Proposition~\ref{prop:LSreduce}, this is deduced from the principal coefficients case.

\begin{proposition}\label{prop:principals=2}
Let $\cA$ be an acyclic cluster variety with principal coefficients.
The spectral sequence $E_{\Filt}$ in mixed Hodge degree $(2,2)$ stabilizes on page $2$. 
The nonzero terms on this page are $E_2^{00, 2} \cong \bigwedge^2 H_{\frozen}$, $E_2^{1(-1), 2} = H_{\GSV} \cong H^0(\Gamma)$ and $E_2^{2(-1), 2} \cong H^1(\Gamma)$.
These terms are in cohomological degrees $2$, $2$, and $3$ respectively. The space  $H^{3,(2,2)}(\cA)$  is spanned by the edge classes.
%We have $H^{\ast,(2,2)}(\cA) = H^{2,(2,2)}(\cA) \oplus H^{3,(2,2)}(\cA)$.  
%\begin{enumerate}
%\item
%The group $H^{2,(2,2)}(\cA)$ decomposes as $\bigwedge^2 H_{\frozen} \oplus H_{\GSV}$.
%\item The group $H^{3,(2,2)}(\cA)$ is naturally isomorphic to $H^1(\Gamma)$.
%\end{enumerate}
\end{proposition}

\begin{proof}
We start by looking at the part of $E_1 = E_{1, \Filt}$ in mixed Hodge degree $(2,2)$. The degree $(2,2)$ terms come from pairs $(D,E)$ with $|D|+|E|=2$, so there are three cases to consider.

If $|D|=2$ and $|E| = 0$, then we are in position $E^{00}$. As we learned from Proposition~\ref{E00}, the part of $E^{00}_r$ in degree $(2,2)$ is $\bigwedge^2 H_{\frozen}$ for every $r$, and the maps from this to every other position in the spectral sequence are $0$.

If $|D|=|E|=1$, then $E \setminus D$ is either a singleton or empty. 
If $E \setminus D$ is singleton, then $\cI(E \setminus D)$ is contractible, so we get no contribution to $E_{1}$. 
If $D=E=\{ i \}$, then $\cI(E \setminus D)$ is empty, so $\tH^{-1}$ is one dimensional. Thus, for each vertex $i$ of $\Gamma$, we get a one dimensional summand of $E^{1(-1)}_{1}$.

Now, suppose that $|D|=0$ and $|E|=2$, say $E = \{ a,b \}$.
If $(a,b)$ is not an edge of $\Gamma$, then $\cI(\{ a,b \})$ is contractible and does not contribute to  $E_{1}$. 
If $(a,b)$ is an edge of $\Gamma$, then $\cI(\{ a,b \})$ is two points, so $\tH^0(\cI(\{ a,b \}))$ is one dimensional. 
Thus, for each edge vertex $\{a,b \}$ of $\Gamma$, we get a one dimensional summand of $E^{2(-1)}_{1}$.

In short,  we have $E^{00,2}_{1} =\bigwedge^2 H_{\frozen}$,  we have $E^{1(-1),2}_{1} =\CC^{\text{Vertices}(\Gamma)}$ and we have  $E^{2(-1),2}_{1} =\CC^{\text{Edges}(\Gamma)}$. 
The map out of position $(0,0)$ will be $0$ on every page, by Proposition~\ref{E00}, so we only need to consider the map between positions $(1,-1)$ and $(2,-1)$; the spectral sequence will then stabilize after that point. 
From Theorem~\ref{thm:E1map}, the map $\CC^{\text{Vertices}(\Gamma)} \longrightarrow \CC^{\text{Edges}(\Gamma)}$ is the map in the chain complex for the simplicial cohomology of $\Gamma$, followed by multiplying the coordinate of edge $(a,b)$ by $B_{ab}$. 
So the cohomology of this map is the simplicial cohomology of $\Gamma$, and we obtain that $E_2^{1(-1), 2} \cong H^0(\Gamma)$ and $E_2^{2(-1), 2} \cong H^1(\Gamma)$.
It is straightforward to identify $H^0(\Gamma)$ with the span of the  GSV classes and $H^1(\Gamma)$ with the span of the edge classes.
\end{proof}

\subsection{The case $s=3$} 
%For $i \in [n]$, let $d_i$ denote the degree of the vertex $i$ in $\Gamma$ and define 
%\[ e_i:= \#\{\text{components of } \Gamma {\setminus} i\} - \#\{\text{components of } \Gamma\}.
%\]
Let $\ell$ denote the number of connected components of $\Gamma$, let $\ell_1$ denote the number of those components that are isolated vertices and let $\Delta(\Gamma)$ denote the number of $3$-cycles in $\Gamma$.
For $i \in [n]$, let $d_i$ denote the degree of vertex $i$ and let $e_i:= \#\{\text{components of } \Gamma {\setminus} i\} - \#\{\text{components of } \Gamma\}$.

\begin{proposition}\label{prop:s=3}
Let $\cA$ be an acyclic cluster variety of really full rank.  Then $H^{\ast,(3,3)}(\cA) = H^{3,(3,3)}(\cA) \oplus H^{4,(3,3)}(\cA)$.  We have $H^{3,(3,3)}(\cA) \cong  \bigwedge^3 H_{\frozen} \oplus H_{\frozen} \wedge H_{\GSV}$, and
\begin{align*}
\dim H^{3,(3,3)}(\cA) &= \binom{m}{3} + m \ell - \ell_1 \\
\dim H^{4,(3,3)}(\cA) &= m \dim H^1(\Gamma) - \sum_i (d_i-e_i+1)   +  \binom{d_i}{2} - \Delta(\Gamma) - \sum_i e_i - \ell_1.
\end{align*}
\end{proposition}
Using Proposition~\ref{prop:LSreduce}, this is deduced from the principal coefficients case.

\begin{proposition}\label{prop:principals=3}
Let $\cA$ be an acyclic cluster variety with principal coefficients.
The spectral sequence $E_{\Filt}$, in mixed Hodge degree $(3,3)$, stabilizes on page $E_3$.
The nonzero terms on this page, and their dimensions, are listed below: 
\[\begin{array}{r@{\ }c@{\ }l@{\quad}r@{\ }c@{\ }l}
 E^{00, 3}_{3, \Filt} &=& \bigwedge^3 H_{\frozen} & \dim  E^{00, 3}_{3, \Filt}  &=& \binom{n}{3}\\ 
 E^{1(-1), 3}_{3, \Filt}  &=& H_{\frozen} \wedge H_{\GSV} & \dim  E^{1(-1), 3}_{3, \Filt}   &=& n \ell - \ell_1 \\ 
 E^{2(-1), 3}_{3, \Filt} &=& H_{\frozen} \wedge H^{3, (2,2)}(\cA) &  \dim  E^{2(-1), 3}_{3, \Filt}  &=& n \dim H^1(\Gamma) - \sum_i (d_i-e_i+1)  \\
 &\cong& \bigoplus_i H^1(\Gamma \setminus \{ i \}) & && \\
\text{and}\ E^{3(-2)}_{3, \Filt} && & \dim  E^{3(-2)}_{3, \Filt} &=& \sum \binom{d_i}{2} - \Delta(\Gamma) - \sum_i e_i - \ell_1.\\
 \end{array} \] 
 These lie in cohomological degrees $3$, $3$, $4$ and $4$ respectively. 
%Their dimensions are $\binom{n}{3}$, $n \ell - \ell_1$, $n \dim H^1(\Gamma) - \sum (d_i-e_i+1)$ and $\sum \tbinom{d_i}{2} - \Delta(\Gamma) - \sum e_i - \ell_1$ respectively
\end{proposition}

These formulas simplify when $\Gamma$ is a tree with at least $2$ vertices; in this case, $\dim E^{00, 3}_{3, \Filt} = \binom{n}{3}$, $\dim E^{1(-1), 3}_{3, \Filt} = n$, $E^{2(-1), 3}_{3, \Filt}=0$ and $\dim E^{3(-2)}_{3, \Filt} = \sum_i \binom{d_i-1}{2}$. 
In particular, if $\Gamma$ is a tree which is not a path, then $H^{4, (3,3)}(\cA)$ is nonzero, so we obtain non-standard cohomology in all of these cases. 

\begin{proof}
We begin by computing $E_1 = E_{1, \Filt}$ in mixed Hodge degree $(3,3)$. We must have $|D|+|E|=3$. 

If $|D|=3$ and $|E|=0$, then we obtain $\bigwedge^3 H_{\frozen}$ in $E^{00,3}_1$.

If $|D|=2$ and $E=|1|$, then $\cI(E \setminus D)$ is either a single point or empty, according to whether $E \cap D = \emptyset$ or $E \subset D$ respectively.
Only the latter case will contribute to $E_1$, and that case contributes in degree $(1,-1)$. 
So $\dim E^{1(-1), 3}_{1} = n(n-1)$, coming from $(D,E)$ of the form $(\{ i, j \}, \{i \})$. 

If $|D|=1$ and $E=|2|$, then the only case where $\cI(E \setminus D)$ is not contractible is when $E = \{ j,k \}$ is an edge of $\Gamma$ and $D = \{ i \}$ for $i \not\in \{ j,k \}$. 
These cases contribute in position $E^{2(-1),3}_1$, coming from $(D,E)$ of the form $(\{ j,k  \}, \{i \})$. 

Finally, we come to the case where $|D|=0$ and $|E|=3$. If $\Gamma_E$ has an isolated vertex, then $\cI(E)$ is contractible, so we only need to consider the cases where $\Gamma_E$ is a $3$-path or a $3$-cycle. These cases give rise to $\CC$ or $\CC^2$ respectively, in $E^{3(-2),3}_1$. The number of embedded $3$-paths in $\Gamma$ is $\sum_i \tbinom{d_i}{2} - 3 \Delta(\Gamma)$, by counting which vertex is the middle of the path. Thus, $\dim E^{3(-2), 3}_{1 } = \sum \tbinom{d_i}{2} - \Delta(\Gamma)$. 

Our next step is to consider maps on the $E^1$ page. The only possible nonzero map is $E^{1(-1), 3}_{1} \longrightarrow E^{2(-1), 3}_{1}$.
The $(\{ i, j \}, \{i \})$ summand on the left is mapped to the sum of the $(\{ j,k  \}, \{i \})$ summands where $k$ ranges over vertices of $\Gamma$ adjacent to $j$, other than $i$. 
Thus,  $E^{1(-1), 3}_{1} \longrightarrow E^{2(-1), 3}_{1}$ splits into a direct sum over fixed values of $D = \{ i \}$, and the cohomology of the $\{ i \}$-summand is $\bigoplus H^0(\Gamma \setminus \{ i \})$ in position $(1,-1)$ and $\bigoplus H^1(\Gamma \setminus \{ i \})$ in position $(2,-1)$. 
We note that the dimensions of these spaces are $\ell n + \sum e_i$ and $n \dim H^1(\Gamma) - \sum (d_i-e_i+1)$ respectively.
In particular, since all maps out of  position $(0,0)$ are the zero map (Proposition~\ref{E00}), we have now computed the dimensions of $E^{00, 3}$ and $E^{2(-1), 3}$ on all pages.

The only possible remaining nonzero map is the map $E^{1(-1), 3}_{2} \longrightarrow E^{3(-2), 3}_{2}$.
Recall that the dimensions of these spaces are $\ell n + \sum e_i$ and $\sum \tbinom{d_i}{2} - \Delta(\Gamma)$.
By the computation of standard cohomology in~\cite{LS},  $H^{3, (3,3)}(\cA) = \bigwedge^3 H_{\frozen} \oplus H_{\frozen} \wedge H_{\GSV}$. 
But we also know that $H^{3, (3,3)}(\cA)$ is filtered with subquotients $E^{00,3}_{3} \cong  \bigwedge^3 H_{\frozen}$ and $E^{1(-1),3}_{3}$.
So $E^{1(-1),3}_{3}$ must be $H_{\frozen} \wedge H_{\GSV}$. 

If $\Gamma'$ is a component of $\Gamma$ which is an isolated vertex $\{ i \}$, then $\dlog y_i \wedge \gamma_{\Gamma'}=0$. 
We claim that $H_{\frozen} \wedge H_{\GSV}$ is the quotient of $H_{\frozen} \otimes H_{\GSV}$ by these relations. 
Indeed, from~\cite{LS}, the standard cohomology of $\cA$ is isomorphic to the ring of differential forms on the big cluster torus generated by the frozen forms and the GSV forms, and it is easy to check that the differential forms $\dlog y_i \wedge \gamma_{\Gamma_j }$, for $\Gamma_j \neq \{ i \}$, are all linearly independent in $\bigwedge^3 \langle \dlog x_i,\ \dlog y_i \rangle$. 
Thus,  $E^{1(-1),3}_{3} \cong H_{\frozen} \wedge H_{\GSV}$ must have dimension $\dim H_{\frozen} \otimes H_{\GSV} - \ell_1 = n\ell - \ell_1$, as claimed, and this space will stabilize on all future pages.

Since $\dim E^{1(-1),2}_{3} =  \ell n + \sum e_i$, we deduce that the rank of the map $E^{1(-1), 3}_{2} \longrightarrow E^{3(-2), 3}_{2}$ is $\sum e_i + \ell_1$. 
This means that the cokernel of this map must have dimension $\sum \tbinom{d_i}{2} - \Delta(\Gamma) - \sum e_i - \ell_1$. 
We have computed the dimension of $E^{3(-2),3}_{3}$, and thus of $E^{3(-2),3}$ on all future pages.
\end{proof}

\section{Examples}\label{sec:examples}
We suppose that $\cA$ has principal coefficients throughout this section.
\subsection{Stars}
Let $\Gamma = Z_n$ denote the star graph with one central vertex $v_0$ joined to $n-1$ vertices of degree $1$ (so $Z_2$ is a single edge).
Every induced subgraph of $Z_n$ is either the empty graph, a positive number of isolated vertices, or the star $Z_m$ for $m \leq n$.
The only nonzero reduced cohomology of $\cI$ in these cases is $\tH^{-1}(\cI(\emptyset)) \cong \CC$ and $\tH^0(\cI(Z_m)) \cong \CC$ whenever $m \geq 2$.
The latter result is because  $\cI(Z_m)$ has the homotopy type of two points (and, in fact, is the disjoint union of an $(m-1)$-simplex and an isolated point).

Since the nonzero  reduced cohomology of $\cI$ lies in $\tH^{-1}$ and $\tH^0$, we see that $E^{ef}_{1, \Filt}$ is only nonzero for $0 \leq e+f \leq 1$. 
Specifically, $E^{e(-e), s}_1$ has a basis indexed by pairs $E \subseteq D$ with $|E| = e$ and $s = |D|+|E|$. 
Meanwhile, $E^{(e)(1-e), s}_1$ has a basis indexed by pairs $(D,E)$ such that $|E \setminus D| \geq 2$, and $v_0 \in E \setminus D$, with $|E| = e$ and $s = |D|+|E|$. 

%By Proposition~\ref{prop:EH}, $\cI(Z_n) \simeq \cI(Z_2) \simeq S^0$ has the homotopy type of two points.  Thus $E_1^{\bullet,-\bullet}$ is the sum of $\tH^{-1}(\emptyset) \simeq \CC$ over all $(B,C)$ such that $C \subset B$ and $E_1^{\bullet,1-\bullet}$ is the sum of $\tH^0(S^0) \simeq \CC$ over all $(B,C)$ such that $C {\setminus} B$ is either an edge or a star $Z_r \subset Z_n$.  These are the only nonzero groups on the $E_1$-page.  It follows that $H^k(\cA) \cong H^{k,(k,k)}(\cA) \oplus H^{k,(k-1,k-1)}(\cA)$.  The only nonzero differentials $\partial_r$ of the $E_r$-page send $E_r^{\bullet,-\bullet}$ to $E_r^{\bullet+r,1-\bullet-r} = E_r^{\bullet,1-\bullet}$.  

It is convenient to encode the dimensions of the spaces in $E_{1}$ by generating functions. 
\[
\begin{array}{lcl}
 \sum_{e,s} \dim E^{e(-e), s}_{1} x^s y^e &=& (1+x+x^2y)^n \\
  \sum_{e,s} \dim E^{e(1-e), s}_{1} x^s y^e &=& xy(1+x)^{n-1} (1+xy)^{n-1} - xy(1+x+x^2y)^{n-1}
 \end{array}
 .\]
 Plugging in $y=1$, we obtain:
 \[
\begin{array}{lcl}
   \sum_{e,s} \dim E^{e(-e), s}_{1} x^s &=& (1+x+x^2)^n \\
   \sum_{e,s} \dim E^{e(1-e), s}_{1} x^s &=& x (1+x)^{2n-2}-x(1+x+x^2)^{n-1} \\
    \end{array}
 .\]

When the spectral sequence eventually stabilizes, the only nonzero cohomology will be in $H^{s, (s,s)}(\cA)$ and $H^{s+1, (s,s)}(\cA)$, with $\dim H^{s, (s,s)}(\cA) = \sum_e \dim E^{e(-e),s}_{\infty}$ and $\dim H^{s+1, (s,s)}(\cA) = \sum_e \dim E^{e(1-e),s}_{\infty}$.
We therefore obtain
\[ \sum_s \dim H^{s, (s,s)}(\cA) x^s - \sum_s \dim H^{s+1, (s,s)}(\cA) x^s = (1+x+x^2)^n - x (1+x)^{2n-2} +x(1+x+x^2)^{n-1}. \]

By Theorem~\ref{thm:LS}, we have 
\[
\sum_{s\geq 0} \dim (H^{s,(s,s)}(\cA)) x^s = (1+x)^{n-1} (1+x+x^2+\cdots+x^{n+1}).
\]
so
\[ \sum_{s \geq 0} \dim H^{s+1,(s,s)}(\cA) x^s = 
x(1+x)^{2n-2}  -  (1+x)^2(1+x+x^2)^{n-1} + (1+x)^{n-1} (1+x+x^2+\cdots+x^{n+1}). \]

%
%Define $g(a,b):=  \sum_{p=0}^{\lfloor b/2 \rfloor} \binom{a}{p} \binom{a-p}{b-2p}$.  Now, the dimensions of $E_1^{\bullet,-\bullet,s}$ are given by
%$$
%\sum_{s \geq 0} \dim(E_1^{\bullet,-\bullet,s}) x^s = (1+x+x^2)^{n}.
%$$ 
%and for $E_1^{\bullet.1-\bullet,s}$ we have
%$$
%\sum_{s \geq 0} \dim(E_1^{\bullet,1-\bullet,s}) x^s = x\left[(1+2x+x^2)^{n-1} - (1+x+x^2)^{n-1}\right].
%$$
%By Theorem~\ref{thm:LS}, we have 
%$$
%\sum_{s\geq 0} \dim (H^{s,(s,s)}(\cA)) x^s = (1+x)^{n-1} (1+x+x^2+\cdots+x^{n+1}).
%$$ 
%Since
%$$
%\dim H^{s+1,(s,s)}(\cA) = \dim E_1^{\bullet,1-\bullet,s} - (\dim E_1^{\bullet,-\bullet,s} - \dim H^{s,(s,s)}(\cA)),
%$$
%we have
%\begin{align*}
%&\sum_{s \geq 0} \dim H^{s+1,(s,s)}(\cA) x^s \\
%&= x(1+2x+x^2)^{n-1}  -  (1+2x+x^2)(1+x+x^2)^{n-1} + (1+x)^{n-1} (1+x+x^2+\cdots+x^{n+1}).
%\end{align*}
For $n = 3$, this generating function $\sum_{s \geq 0} \dim H^{s+1,(s,s)}(\cA) x^s$ is 0.  For $n =4$, we get $x^3+2x^4+x^5$, and for $n = 5$ we get $3x^3 + 11 x^4 + 16 x^5 + 11 x^6 + 3 x^7$.  The palindromicity of these polynomials agrees with the curious Lefschetz theorem \cite{LS}.  The coefficient of $x^3$ is equal to $\binom{n-2}{2}$, agreeing with Proposition~\ref{prop:principals=3}.

\subsection{Paths}  Let $\Gamma = A_n$ be a path on $n$ vertices.  By Theorem~\ref{thm:LS}(4), we have $H^*(\cA) = \bigoplus_s H^{s,(s,s)}(\cA)$.  The spectral sequence $E_{r,\Filt}^{pq}$ of the filtered complex $F^\bullet G^\bullet$ is however quite complicated.  The $E_1$-page is the direct sum $E_1^{\bullet, -\bullet} \oplus E_1^{\bullet, 1-\bullet} \oplus E_1^{\bullet, 2-\bullet} \oplus \cdots$.  The group $E_1^{\bullet, -\bullet}$ is the sum of $\tH^{-1}(\emptyset) \simeq \CC$ over all $(B,C)$ such that $C \subset B$.  The group $E_1^{\bullet,t-\bullet}$ is the sum of $\tH^{t-1}(S^{t-1}) \simeq \CC$ over all $(B,C)$ such that $\cI(\Gamma_{C {\setminus} B})$ is a sphere of dimension $t-1$.  By Proposition~\ref{prop:join} and Corollary~\ref{cor:path}, this happens exactly when $\Gamma_{C {\setminus} B}$ is a disjoint union of paths $P_{a_1},P_{a_2},\ldots,P_{a_r}$ where $a_i \in \{3k_i+1,3k_i +2 \}$ and $k_1+k_2+\cdots+k_r + r = t$.
It is not clear to us on what page $E_{\Filt}$ collapses for the path.

\subsection{Forests}
Let $\Gamma$ be a forest, so every induced subgraph $F$ of $\Gamma$ is also a forest. By Proposition~\ref{prop:EH}, $\cI(F)$ is either contractible, or a sphere.
So, for all subsets $D$ and $E$ of $[n]$, $\cI(E \setminus D)$ only has reduced cohomology in at most one degree, and this reduced cohomology has dimension $1$.
So the $E_1$ page has one basis element indexed by each pair $(D,E)$ for which $\cI(E \setminus D)$ is a sphere.

We can also understand the maps on the $E_1$-page, using Theorem~\ref{thm:E1map}.
Let $(x,y)$ be an edge of $\Gamma$, let $F$ be a subforest of $\Gamma \setminus \{ x,y \}$ and let $F'$ be the forest induced on the vertices of $F \cup \{ x,y \}$. 
Let $T'$ be the component of $F'$ containing the edge $(x,y)$ and let $T$ be the forest $T' \setminus \{ x,y \}$.
We'd like to understand the Mayer-Vietores map $\widetilde{H}^r(\cI(F)) \to \widetilde{H}^{r+1}(\cI(F'))$; by Proposition~\ref{prop:join}, this is the join of the map $\widetilde{H}^r(\cI(T)) \to \widetilde{H}^{r+1}(\cI(T'))$ with a fixed other factor coming from the other components of $F'$.

One can check that all of the homotopy equivalences in Proposition~\ref{prop:EH} commute with Mayer-Vietores maps. 
Repeatedly using these equivalences to eliminate vertices of $F'$ other than $x$ and $y$, we may assume that $F'$ is a path of length $a+b+1$, with $(x,y)$ the $(a+1)$-st edge, and $0 \leq a,b \leq 2$. 
Checking these cases, we obtain:
\begin{prop} \label{forest cases}
Let $T'$ be a path of length $a+b+1$ and let $(x,y)$ be the $(a+1)$-st edge of $T'$, for $0 \leq a,b \leq 2$. Let $T$ be the forest $T' \setminus \{ x,y \}$. Then the Mayer-Vietores map $\widetilde{H}^{\ast}(\cI(T)) \to \widetilde{H}^{\ast+1}(\cI(T'))$ is an isomorphism $\widetilde{H}^{-1} \to \widetilde{H}^0$ if $a=b=0$, and is the $0$ map otherwise. 
\end{prop}

\begin{proof}
If $a$ or $b=1$, then $T$ has an isolated vertex, so $\cI(T)$ is contractible by Corollary~\ref{cor:isolated}. If $(a,b) = (0,2)$ or $(2,0)$, then $\cI(T')$ is contractible. The remaining cases are $(a,b) = (0,0)$ or $(2,2)$.

When $(a,b) = (0,0)$, the forest $T$ is the empty set so $\cI(T) = \emptyset$ and $\widetilde{H}^{-1}(\cI(T)) = \CC$; meanwhile, the forest $T'$ is a single edge, so $\cI(T')$ is two points and $\widetilde{H}^0(\cI(T')) = \CC$. The Mayer-Vietores cover in question is covering the two points $\{ x, y \}$ of  $\cI(T')$  by $\{x \}$ and by $\{y \}$, and the map is an isomorphism.

When $(a,b) = (2,2)$, the only non-vanishing cohomology groups are $\widetilde{H}^1(\cI(T)) \cong \widetilde{H}^1(\cI(T'))$, so these cannot be connected by a Mayer-Vietores map.
\end{proof}

This makes it tractable to write down $E^1_{\Filt}$ in practice.
We have little understanding of the maps on the later pages.


\begin{thebibliography}{xxx}
\bibitem[Ara]{Ara} D.~Arapura.  The Leray spectral sequence is motivic.
Invent. Math. 160 (2005), no. 3, 567--589.
\bibitem[ABL]{ABL} N.~Arkani-Hamed, Y.~Bai, and T.~Lam. Positive geometries and canonical forms, JHEP 2017, Article number: 39 (2017).
\bibitem[ABCGPT]{book} N.~Arkani-Hamed, J.~Bourjaily, F.~Cachazo, A.~Goncharov, A.~Postnikov, and J.~Trnka. Grassmannian geometry of scattering amplitudes. Cambridge University Press, Cambridge, 2016. ix+194 pp.
\bibitem[AHL]{AHL} N.~Arkani-Hamed, S.~He, and T.~Lam. Cluster configuration spaces of finite type, SIGMA 17 (2021), 092, 41 pages.
\bibitem[BFZ]{CA3} A. Berenstein, S.~Fomin, and A.~Zelevinsky. Cluster algebras. III. Upper bounds and double Bruhat cells. Duke Math. J. 126 (2005), no. 1, 1--52.
\bibitem[Del]{Del} P. Deligne. Th\'eorie de Hodge. I, Actes du Congr\`es International des Math\'ematiciens, Gauthier-Villars, Paris, 1971, pp. 425--430.
\bibitem[EH]{EH} R.~Ehrenborg and G.~Hetyei. The topology of the independence complex. European J. Combin. 27 (2006), no. 6, 906--923.
\bibitem[FZ]{CA1} S.~Fomin and A.~Zelevinsky, Cluster algebras. I. Foundations, JAMS (2002), no. 2, 497--529.
\bibitem[GL]{GL} P.~Galashin and T.~Lam, Positroids, knots, and $q,t$-Catalan numbers, preprint, 2020; {\tt arXiv:2012.09745}.
\bibitem[GSV]{GSV} M.~Gekhtman, M.~Shapiro, and A.~Vainshtein. Cluster algebras and Poisson geometry. Mathematical Surveys and Monographs, 167. American Mathematical Society, Providence, RI, 2010. xvi+246 pp.
\bibitem[Hag]{Hag} J.~Haglund.  The $q,t$-Catalan numbers and the space of diagonal harmonics. With an appendix on the combinatorics of Macdonald polynomials. University Lecture Series, 41. American Mathematical Society, Providence, RI, 2008. viii+167 pp.
\bibitem[HLZ]{HLZ}  A. Huang, B. H. Lian, and X. Zhu. Period integrals and the Riemann-Hilbert correspondence, J. Differential Geom. 104 (2016), no. 2, 325--369.
\bibitem[Koz]{Koz} D. N.~Kozlov, Complexes of directed trees. J. Combin. Theory Ser. A 88 (1999), no. 1, 112--122.
\bibitem[LS]{LS} T.~Lam and D.~Speyer.  Cohomology of cluster varieties. I. Locally acyclic case.  Algebra Number Theory, to appear.
\bibitem[LT]{LT} T.~Lam and N.~Templier. 
The mirror conjecture for minuscule flag varieties, preprint, 2017; {\tt arXiv:1705.00758}.
\bibitem[Mul]{Mul} G.~Muller. Locally acyclic cluster algebras. Adv. Math. 233 (2013), 207--247. 
\bibitem[MS]{MS} G.~Muller and D.~E.~Speyer. Cluster algebras of Grassmannians are locally acyclic. Proc. Amer. Math. Soc. 144 (2016), no. 8, 3267--3281.
\bibitem[MT1]{MTarXiv} M.~Marietti and D.~Testa.  Cores of simplicial complexes, preprint 2007; {\tt arXiv:math/0703351}.
\bibitem[MT2]{MTpub} M.~Marietti and D.~Testa.  Cores of simplicial complexes. Discrete Comput. Geom. 40 (2008), no. 3, 444--468.
\bibitem[Pet]{Pet} D.~Petersen. A spectral sequence for stratified spaces and configuration spaces of points. Geom. Topol. 21 (2017), no. 4, 2527--2555.
%\bibitem[PS]{PS}  C.~Peters and J.~Steenbrink, Mixed Hodge structures. Ergebnisse der Mathematik und ihrer
Grenzgebiete. 3. Folge., vol. 52, Springer-Verlag, Berlin, 2008.
\bibitem[Spe]{SpeyerFinGenCounterExample} D.~E.~Speyer, An infinitely generated upper cluster algebra, preprint 2013, {\tt arXiv:1305.6867}.
\bibitem[Stacks]{stacks-project} The Stacks project, 2021, 
  \url{https://stacks.math.columbia.edu}.\end{thebibliography}
\end{document}